\documentclass[12pt,leqno]{amsart}

\usepackage{amscd}
\usepackage{amsmath}
\usepackage{amsfonts}
\usepackage{amssymb}
\usepackage{url}
\usepackage[small,nohug,heads=LaTeX]{diagrams}
\diagramstyle[labelstyle=\scriptstyle]

\title{Noether's problem for central extensions of metacyclic $p$-groups}
\author{Ivo M. Michailov, Ivan S. Ivanov}
\address{Faculty of Mathematics and Informatics, Constantin Preslavski University, Universitetska str. 115, 9700 Shumen, Bulgaria}
\email{ivo\_michailov@yahoo.com;slaveikov@abv.bg}
\date{\today}
\thanks{This work is partially supported by a project No RD-05-156/25.02.2011 of Shumen University}
\begin{document}
\baselineskip 20pt
\begin{abstract}
Let $K$ be a field and $G$ be a finite group. Let $G$ act on the
rational function field $K(x(g):g\in G)$ by $K$ automorphisms
defined by $g\cdot x(h)=x(gh)$ for any $g,h\in G$. Denote by $K(G)$
the fixed field $K(x(g):g\in G)^G$. Noether's problem then asks
whether $K(G)$ is rational over $K$. In [M. Kang, Noether's problem
for metacyclic $p$-groups, Adv. Math. 203(2005), 554-567], Kang
proves the rationality of $K(G)$ over $K$ if $G$ is any metacyclic
$p$-group and $K$ is any field containing enough roots of unity. In
this paper, we give a positive answer to the Noether's problem for
all central group extensions of the general metacyclic $p$-group,
provided that $K$ is infinite and it contains sufficient roots of
unity.
\end{abstract}

\maketitle
\newcommand{\Gal}{{\rm Gal}}
\newcommand{\Ker}{{\rm Ker}}
\newcommand{\GL}{{\rm GL}}
\newcommand{\Br}{{\rm Br}}
\newcommand{\lcm}{{\rm lcm}}
\newcommand{\ord}{{\rm ord}}
\numberwithin{equation}{section}

\section{Introduction}
\label{1}

Let $K$ be a field and $G$ be a finite group. Let $G$ act on the
rational function field $K(x(g):g\in G)$ by $K$ automorphisms
defined by $g\cdot x(h)=x(gh)$ for any $g,h\in G$. Denote by $K(G)$
the fixed field $K(x(g):g\in G)^G$. {\it Noether's problem} then
asks whether $K(G)$ is rational ($=$ purely transcendental) over
$K$. It is related to the inverse Galois problem, to the existence
of generic $G$-Galois extensions over $k$, and to the existence of
versal $G$-torsors over $k$-rational field extensions \cite[33.1,
p.86]{Sw,Sa1,GMS}.

The following well-known theorem gives a positive answer to the
Noether's problem for abelian groups.

\newtheorem{t1.1}{Theorem}[section]
\begin{t1.1}\label{t1.1}
{\rm (Fischer }\cite[Theorem 6.1]{Sw}{\rm )} Let $G$ be a finite
abelian group of exponent $e$. Assume that {\rm (i)} either char $K
= 0$ or char $K > 0$ with char $K\nmid e$, and {\rm (ii)} $K$
contains a primitive $e$-th root of unity. Then $K(G)$ is rational
over $K$.
\end{t1.1}

Swan's paper \cite{Sw} also gives a survey of many results related
to the Noether's problem for abelian groups. In the same time, just
a handful of results about Noether's problem are obtained when the
groups are nonabelian. The reader is referred to
\cite{CK,Ka1,HuK,Ka2} for previous results of Noether's problem for
$p$-groups.

We state now the result obtained recently by Kang \cite{Ka1} about
the Noether's problem for metacyclic $p$-groups:

\newtheorem{t1.2}[t1.1]{Theorem}
\begin{t1.2}\label{t1.2}
{\rm (Kang }\cite[Theorem 1.5]{Ka1}{\rm )} Let $G$ be a metacyclic
$p$-group with exponent $p^e$, and let $K$ be any field such that
{\rm (i)} char $K = p$, or {\rm (ii)} char $K \ne p$ and $K$
contains a primitive $p^e$-th root of unity. Then $K(G)$ is rational
over $K$.
\end{t1.2}

It is still an open question whether the above result could be
extended for all series of $2$-generator $p$-groups or meta-abelian
groups, that have quotient groups isomorphic to the metacyclic
$p$-groups. However, we should not ''over-generalize'' Theorem
\ref{t1.2}, because Saltman proves the following result.

\newtheorem{t1.3}[t1.1]{Theorem}
\begin{t1.3}\label{t1.3}
{\rm (Saltman }\cite{Sa2}{\rm )} For any prime number $p$ and for
any field $K$ with char $K \ne p$ (in particular, $K$ may be an
algebraically closed field), there is a meta-abelian $p$-group $G$
of order $p^9$ such that $K(G)$ is not rational over $K$.
\end{t1.3}

Among the known results of Noether's problem for non-abelian
$p$-groups, assumptions on the existence of ``enough" roots of unity
always arose. In fact, even when $G$ is a non-abelian $p$-group of
order $p^3$ where $p$ is an odd prime number, it is not known how to
find a necessary and sufficient condition such as $\mathbb Q(G)$ is
rational over $\mathbb Q$ (see \cite{Ka3}). Thus it will be
desirable if we can weaken the assumptions on the existence of roots
of unity.

The purpose of this paper is to extend Theorem \ref{t1.2} for all
central $p$-extensions of the general metacyclic $p$-group. However,
some additional assumptions will appear in the statements of our
results so that we guarantee the existence of the groups we are
going to consider.

Let $G$ be any metacyclic $p$-group generated by two elements
$\sigma$ and $\tau$ with relations
$\sigma^{p^a}=1,\tau^{p^b}=\sigma^{p^c}$ and
$\tau^{-1}\sigma\tau=\sigma^{\varepsilon+\delta p^r}$ where
$\varepsilon=1$ if $p$ is odd, $\varepsilon=\pm 1$ if $p=2$,
$\delta=0,1$ and $a,b,c,r\geq 0$ are subject to some restrictions.
For the the description of these restrictions see e.g. \cite[p.
564]{Ka1}. Note that if $\delta=0$ and $\varepsilon=1$, then $G$ is
abelian group generated by two elements.

Our first main result is the following.

\newtheorem{t1.4}[t1.1]{Theorem}
\begin{t1.4}\label{t1.4}
Let $G$ be an abelian metacyclic $p$-group generated by two elements
$\sigma$ and $\tau$ with relations $\sigma^{p^a}=1,\tau^{p^b}=1$ and
$\tau^{-1}\sigma\tau=\sigma$. Assume that $\widetilde G$ is a
central extension of $G$, i.e., we have the following group
extension
\begin{equation*}
1\longrightarrow C\longrightarrow \widetilde G\longrightarrow G\cong
C_{p^a}\times C_{p^b} \longrightarrow 1,
\end{equation*}
where $C\leq Z(\widetilde G)$. Let $p^t$ be the exponent of $C$, let
$a\geq b\geq t$ and let the pre-image of
$[\sigma,\tau]=\sigma^{-1}\tau^{-1}\sigma\tau$ in $\widetilde G$ is
of order $p^t$. Let $e=\max\{a,2t\}$. Assume that {\rm (i)} $char
K=p$ or {\rm (ii)} $char K\ne p$, $K$ is infinite, and $K$ contains
a primitive $p^e$-th root of unity. Then $K(\widetilde G)$ is
rational over $K$.
\end{t1.4}

Next, we consider the case when $G$ is a nonabelian metacyclic
$p$-group. The second main result of this paper is the following
theorem that concerns the central extensions of $G$.

\newtheorem{t1.5}[t1.1]{Theorem}
\begin{t1.5}\label{t1.5}
Let $G$ be a nonabelian metacyclic $p$-group generated by two
elements $\sigma$ and $\tau$ with relations
$\sigma^{p^a}=1,\tau^{p^b}=\sigma^{p^c}$ and
$\tau^{-1}\sigma\tau=\sigma^{\varepsilon+p^r}$, where $1\leq c\leq
a,r\leq\min\{b,c\}$; $\varepsilon=1$ if $p$ is odd and
$\varepsilon=\pm 1$ if $p=2$. Assume that $\widetilde G$ is a
central extension of the group $G$ by the cyclic group $C_{p^t}$,
i.e., we have the following group extension
\begin{equation*}
1\longrightarrow C_{p^t}\longrightarrow \widetilde G\longrightarrow
G \longrightarrow 1,
\end{equation*}
where $C_{p^t}\leq Z(\widetilde G)$. Let $a\geq t, b\geq t$ and let
the pre-image of
$\sigma^{-(k-1)}[\sigma,\tau]=\sigma^{-k}\tau^{-1}\sigma\tau$ in
$\widetilde G$ is of order $p^t$. Let $p^m=\exp(\widetilde G)$ and
$e=\max\{m,r+t\}$. Assume that {\rm (i)} $char K=p$ or {\rm (ii)}
$char K\ne p$, $K$ is infinite, and $K$ contains a primitive
$p^e$-th root of unity. Then $K(\widetilde G)$ is rational over $K$.
\end{t1.5}

We can generalize Theorem \ref{t1.5} for all central extensions
(i.e., not necessarily cyclic extensions), with the expense of a
stronger requirement for the root of unity. One can see from our
proof that this new condition can be weakened in most of the cases
we consider. However, we prefer for simplicity to refrain from
giving these details in the statement of the following result.

\newtheorem{c1.6}[t1.1]{Corollary}
\begin{c1.6}\label{c1.6}
Let $G$ be a nonabelian metacyclic $p$-group generated by two
elements $\sigma$ and $\tau$ with relations
$\sigma^{p^a}=1,\tau^{p^b}=\sigma^{p^c}$ and
$\tau^{-1}\sigma\tau=\sigma^{\varepsilon+p^r}$, where
$\varepsilon=1$ if $p$ is odd and $\varepsilon=\pm 1$ if $p=2$.
Assume that $\widetilde G$ is a central extension of the group $G$,
i.e., we have the following group extension
\begin{equation*}
1\longrightarrow C\longrightarrow \widetilde G\longrightarrow G
\longrightarrow 1,
\end{equation*}
where $C\leq Z(\widetilde G)$. Let $p^t$ be the exponent of $C$, let
$a\geq t, b\geq t$ and let the pre-image of
$\sigma^{-(k-1)}[\sigma,\tau]=\sigma^{-k}\tau^{-1}\sigma\tau$ in
$\widetilde G$ is of order $p^t$. Let $e=a+b+t-c$. Assume that {\rm
(i)} $char K=p$ or {\rm (ii)} $char K\ne p$, $K$ is infinite, and
$K$ contains a primitive $p^e$-th root of unity. Then $K(\widetilde
G)$ is rational over $K$.
\end{c1.6}

If $r\geq t$, we are able to weaken the condition for the roots of
unity in Theorem \ref{t1.5}.

\newtheorem{t1.7}[t1.1]{Theorem}
\begin{t1.7}\label{t1.7}
Let $G$ be a nonabelian metacyclic $p$-group generated by two
elements $\sigma$ and $\tau$ with relations
$\sigma^{p^a}=1,\tau^{p^b}=\sigma^{p^c}$ and
$\tau^{-1}\sigma\tau=\sigma^{\varepsilon+p^r}$, where $1\leq c\leq
a,r\leq\min\{b,c\}$; $\varepsilon=1$ if $p$ is odd and
$\varepsilon=\pm 1$ if $p=2$. Assume that $\widetilde G$ is a
central extension of the group $G$ by the cyclic group $C_{p^t}$,
i.e., we have the following group extension
\begin{equation*}
1\longrightarrow C_{p^t}\longrightarrow \widetilde G\longrightarrow
G \longrightarrow 1,
\end{equation*}
where $C_{p^t}\leq Z(\widetilde G)$. Let $t\leq r$ and let the
pre-image of
$\sigma^{-(k-1)}[\sigma,\tau]=\sigma^{-k}\tau^{-1}\sigma\tau$ in
$\widetilde G$ is of order $p^t$. Let $p^m=\exp(\widetilde G)$.
Assume that {\rm (i)} $char K=p$ or {\rm (ii)} $char K\ne p$, $K$ is
infinite, and $K$ contains a primitive $p^m$-th root of unity. Then
$K(\widetilde G)$ is rational over $K$.
\end{t1.7}

We organize this paper as follows. In Section \ref{2} we recall some
preliminaries which will be used in the proofs of Theorems
\ref{t1.4}, \ref{t1.5}, \ref{t1.7} and Corollary \ref{c1.6}. We give
the presentations of all cyclic central $p$-extensions of the
metacyclic $p$-groups in Section \ref{3}. We also prove some key
results in Sections \ref{2} and \ref{3} that will aid our
investigations later in our work. One of these results is Theorem
\ref{t2.8} which is of interest itself. The proofs of Theorems
\ref{t1.4}, \ref{t1.5} and \ref{t1.7} are given, respectively, in
Sections \ref{4}, \ref{5} and \ref{7}. In Section \ref{6} the proof
of Corollary \ref{c1.6} is given.

\medskip
{\bf Standing notations.} A field extension $L$ of $K$ is rational
over $K$ if $L$ is purely transcendental over $K$. Recall that
$K(G)$ denotes $K(x(g):g\in G)^G$ where $g\cdot x(h)=x(gh)$ for any
$g,h\in G$. A group $G$ is called metacyclic, if $G$ can be
generated by two elements $\sigma$ and $\tau$, and one of them
generates a normal subgroup of $G$. The exponent of a finite group
$G$ is $\lcm\{\ord(g):g\in G\}$ where $\ord(g)$ is the order of $g$.
Two extension fields $L_1$ and $L_2$ of $K$ with $G$-actions are
$G$-isomorphic if there is an isomorphism $\varphi:L_1\to L_2$ over
$K$ such that $\varphi(\sigma\cdot u)=\sigma\cdot\varphi(u)$ for any
$\sigma\in G$, any $u\in L_1$.

\section{Generalities}
\label{2}

We list several results which will be used in the sequel.

\newtheorem{t2.1}{Theorem}[section]
\begin{t2.1}\label{t2.1}
{\rm (}\cite[Theorem 1]{HK}{\rm )} Let $G$ be a finite group acting
on $L(x_1,\dots,x_m)$, the rational function field of $m$ variables
over a field $L$ such that
\begin{description}
    \item [(i)] for any $\sigma\in G, \sigma(L)\subset L;$
    \item [(ii)] the restriction of the action of $G$ to $L$ is
    faithful;
    \item [(iii)] for any $\sigma\in G$,
    \begin{equation*}
\begin{pmatrix}
\sigma(x_1)\\
\vdots\\
\sigma(x_m)\\
\end{pmatrix}
=A(\sigma)\begin{pmatrix}
x_1\\
\vdots\\
x_m\\
\end{pmatrix}
+B(\sigma)
\end{equation*}
where $A(\sigma)\in\GL_m(L)$ and $B(\sigma)$ is $m\times 1$ matrix
over $L$. Then there exist $z_1,\dots,z_m\in L(x_1,\dots,x_m)$ so
that $L(x_1,\dots,x_m)^G=L^G(z_1,\dots,z_m)$ and $\sigma(z_i)=z_i$
for any $\sigma\in G$, any $1\leq i\leq m$.
\end{description}
\end{t2.1}

\newtheorem{t2.2}[t2.1]{Theorem}
\begin{t2.2}\label{t2.2}
{\rm (}\cite[Theorem 3.1]{AHK}{\rm )} Let $G$ be a finite group
acting on $L(x)$, the rational function field of one variable over a
field $L$. Assume that, for any $\sigma\in G,\sigma(L)\subset L$ and
$\sigma(x)=a_\sigma x+b_\sigma$ for any $a_\sigma,b_\sigma\in L$
with $a_\sigma\ne 0$. Then $L(x)^G=L^G(z)$ for some $z\in L[x]$.
\end{t2.2}

\newtheorem{t2.3}[t2.1]{Theorem}
\begin{t2.3}\label{t2.3}
{\rm (}\cite[Theorem 1.7]{CK}{\rm )} If $char K=p>0$ and $\widetilde
G$ is a finite $p$-group, then $K(G)$ is rational over $K$.
\end{t2.3}

Let $\Br(K)$ denote the Brauer group of a field $K$, and $\Br_N(K)$
its $N$-torsion subgroup for any $N>1$. Following Roquette
\cite{Ro}, if $\gamma=[B]\in\Br(K)$ is the class of a $K$-central
simple algebra $B$ and $m\geq 1$ is a multiple of the index of $B$,
then $F_m(\gamma)$ denotes the $m$-th Brauer field of $\gamma$.
Moreover, $F_m(\gamma)/K$ is a regular extension of transcedence
degree $m-1$, which is rational if and only if $\gamma$ is trivial.
The following Theorem was essentially obtained by Saltman in
\cite[p. 541]{Sa3} and proved in detail by B. Plans \cite[Prop.
7]{Pl}.

\newtheorem{t2.4}[t2.1]{Theorem}
\begin{t2.4}\label{t2.4}
{\rm (Saltman }\cite[Proposition 7]{Pl}{\rm )} Let $1\rightarrow
C\rightarrow H\rightarrow G\rightarrow 1$ be a central extension of
finite groups, representing an element $\varepsilon\in H^2(G,C)$.
Let $K$ be an infinite field and let $N$ denote the exponent of $C$.
Assume that $N$ is prime to the characteristic of $K$ and that $K$
contains $\mu_N$ -- the group of $N$-th roots of unity. Let be given
a decomposition $C\cong \mu_{N_1}\times\cdots\times \mu_{N_r}$, and
let the corresponding isomorphism $H^2(G,C)\cong
\oplus_iH^2(G,\mu_{N_i})$ map $\varepsilon$ to $(\varepsilon)_i$.
Let also be given a faithful subrepresentation $V$ of the regular
representation of $G$ over $K$, and let
$\gamma_i\in\Br_N(K(V)^G)\subset \Br(K(V)^G)$ be the inflation of
$\varepsilon_i$ with respect to the isomorphism
$G\cong\Gal(K(V)/K(V)^G)$. Then
\begin{equation*}
K(H) \text{ is rational over the } K(V)^G- \text{ free compositum }
F_m(\gamma_1)\cdots F_m(\gamma_r),
\end{equation*}
where $m$ denotes the order of $G$.
\end{t2.4}

The following key result was suggested by Plans \cite[Proposition
9]{Pl}, but the proof of case 9a, which we need in our paper, is
somewhat hard to extract. We will reprove this result for $p$-group
extensions, following a different approach that will allow us to
weaken the condition for the roots of unity.

\newtheorem{t2.8}[t2.1]{Theorem}
\begin{t2.8}\label{t2.8}
{\rm (Plans }\cite[Proposition 9a]{Pl}{\rm )} Let $1\rightarrow
C\rightarrow H\rightarrow G\rightarrow 1$ be a central extension of
finite $p$-groups for a prime $p$. Let $G'$ (resp. $H'$) be the
derived subgroup of $G$ (resp. $H$), and let $p^n$ (resp. $p^e$)
denote the exponent of $C$ (resp. $G/G'$). Assume that $H'\cap
C=\{1\}$.

{\rm (i)} Let $C=\mu_{p^n}$ and set
$m=\max\{i:\mu_{p^i}\cap\mu_{p^n}\ne\{1\},\mu_{p^i}\leq H/H'\}$. Let
$K$ be an infinite field containing the $p^m$-th roots of unity.
Then $K(H)$ is rational over $K(G)$.

{\rm (ii)} Let $K$ be an infinite field containing the $p^{n+e}$-th
roots of unity. Then $K(H)$ is rational over $K(G)$.
\end{t2.8}
\begin{proof}
Let $\varepsilon\in H^2(G,C)$ be the class of $1\rightarrow
C\rightarrow H\rightarrow G\rightarrow 1$. Define $k=K(G)=K(W)^G$,
where $W$ denotes the regular representation of $G$ over $K$. For
each $\mu_j\subset K^*$ we have an inflation map
$\inf:H^2(G,\mu_j)\to\ _j\Br(k)\subset\Br(k)$ coming from
$\Gal(K(W)/k)\cong G$. Let us consider a decomposition
$C\cong\mu_{p^{s_1}}\times\cdots\times\mu_{p^{s_r}}$, and let
$\varepsilon$ map to $(\varepsilon_i)_i$ via the isomorphism
$H^2(G,C)\cong \bigoplus_iH^2(G,\mu_{p^{s_i}})$. Take
$\gamma_i=\inf(\varepsilon_i)$.

(i) Let $C=\mu_{p^n}$. From $H'\cap C=\{1\}$ it follows that $H'$ is
isomorphic to $G'$. Therefore, we have the following commutative
diagram with exact rows:
\begin{diagram}
1&\rTo&\mu_{p^n}&\rTo&H&\rTo&G&\rTo&1\\
 &   &||&    &\dTo        &    &\dTo\\
1&\rTo&\mu_{p^n}&\rTo&H/H'&\rTo&G/G'&\rTo&1.
\end{diagram}
Let $\varepsilon'\in H^2(G/G',\mu_{p^n})$ be the class of the second
row of the above diagram. Then we have the homomorphism
$\delta:H^2(G/G',\mu_{p^n})\to H^2(G,\mu_{p^n})$ such that
$\delta(\varepsilon')=\varepsilon$. Note that $\varepsilon'$ is not
always split, since $\mu_{p^n}$ may be contained in a bigger cyclic
subgroup of $H/H'$, and we know that extensions of that kind are not
split. Hence $\varepsilon$ is not always split and we can not apply
Theorem \ref{t2.4} at this point.

Thus, we need to take a roundabout approach, involving the theory of
obstructions of Galois embedding problems. Let us recall some
basics. The embedding problem related to $K(W)/k$ and $\varepsilon$
consists in determining whether there exists a Galois algebra (or a
Galois extension) $L$, such that $K(W)$ is contained in $L$, $H$ is
isomorphic to $\Gal(L/k)$, and the homomorphism of restriction of
$L$ on $K(W)$ coincides with the epimorphism $H\to G$. The element
$\gamma=\inf(\varepsilon)\in\Br(k)$ is called the obstruction to
this embedding problem. It is known that the splitting of the
obstruction in the Brauer group is necessary and sufficient
condition for the solvability of such types of embedding problems
(called Brauer problems, see \cite{Le}).

Now, we have that $\inf(\varepsilon')$ is the obstruction to the
embedding problem given by $\varepsilon'$. Let us write the
decomposition of $H/H'$ as a direct product of cyclic groups:
$H/H'\cong\mu_{p^{m_1}}\times\cdots\times\mu_{p^{m_s}}$. For any
$m_j$, set $\mu_{p^{i_j}}=\mu_{p^{m_j}}\cap\mu_{p^n}$, and consider
the embedding problem given by the group extension $1\rightarrow
\mu_{p^{i_j}}\rightarrow \mu_{p^{m_j}}\rightarrow
\mu_{p^{m_j-i_j}}\rightarrow 1$ and some $\mu_{p^{m_j-i_j}}$
extension $L/k$.

Since $k$ contains a primitive $p^{m_j-i_j}$-th root of unity, we
may assume that $L/k=k(\root p^{m_j-i_j}\of{a})/k$ for a certain
element $a\in k^*\setminus k^{*p}$. In this way, we obtain an
equivalent embedding problem given by $1\rightarrow
\mu_{p^{m_j-1}}\rightarrow \mu_{p^{m_j}}\rightarrow \mu_p\rightarrow
1$ and $k(\root p\of{a})/k$. It is well-known (see \cite[Theorem
11]{Al}) that the obstruction to the latter problem is equal to the
Hilbert symbol $(a,\zeta_{p^{m_j-1}})_{p^{m_j-1}}$ which is trivial
for any $a$, since $\zeta_{p^{m_j}}=\root p\of{\zeta_{p^{m_j-1}}}$
is in $k$ (recall that $m=\max\{m_j:1\leq j\leq s\}$). Therefore,
there exists a solution to the embedding problem given by
$\varepsilon'$ that is isomorphic to the $k$-compositum of the
solutions to the embedding problems given by $1\rightarrow
\mu_{p^{i_j}}\rightarrow \mu_{p^{m_j}}\rightarrow
\mu_{p^{m_j-i_j}}\rightarrow 1$ for $1\leq j\leq s$. Hence
$\inf(\varepsilon')=0$.

Next, observe that the inflation map $\inf:
H^2(G/G',\mu_{p^n})\to\Br(k)$ factors through
$\delta:H^2(G/G',\mu_{p^n})\to H^2(G,\mu_{p^n})$, i.e., we have the
following commutative diagram
\begin{diagram}
   &H^2(G/G',\mu_{p^n})&\rTo^{\delta}&H^2(G,\mu_{p^n})&\qquad\qquad\qquad\qquad\qquad\qquad\qquad\\
\qquad\qquad\qquad\qquad\qquad\qquad&\dTo_{\inf}&&\dTo_{\inf}&\\
&\Br(k)&\rTo^{i}&\Br(k),&
\end{diagram}
where $i$ is an injection. Therefore,
$\gamma=\inf(\varepsilon)=\inf\delta(\varepsilon')=\inf(\varepsilon')=0$.
From Theorem \ref{t2.4} now it follows that $K(H)$ is rational over
$K(G)$.

(ii) Let us consider a decomposition
$C\cong\mu_{p^{s_1}}\times\cdots\times\mu_{p^{s_r}}$. If
$\mu_{p^{s_i}}\cap\mu_{p^{m_i}}\ne\{1\}$ for some cyclic subgroup
$\mu_{p^{m_i}}$ of $H/H'$, then $m_i\leq n+e$, as $p^e$ is the
exponent of $G/G'$. Since $K$ contains a primitive $p^{n+e}$-th root
of unity, we can apply case (i) and get that all the $\gamma_i$'s
are trivial. Therefore, $K(H)$ is rational over $K(G)$ by Theorem
\ref{t2.4}.
\end{proof}

\newtheorem{t2.6}[t2.1]{Theorem}
\begin{t2.6}\label{t2.6}
{\rm (Kang }\cite[Theorem 4.1]{Ka1}{\rm )} Let $p$ be a prime
number, $m,n$ and $r$ are positive integers, $k=1+p^r$ if $(p,r)\ne
(2,1)$ (resp. $k=-1+2^r$ with $r\geq 2$). Let $G$ be a split
metacyclic $p$-group of order $p^{m+n}$ and exponent $p^e$ defined
by $G=\langle\sigma,\tau:
\sigma^{p^m}=\tau^{p^n}=1,\tau^{-1}\sigma\tau=\sigma^k\rangle$. Let
$K$ be any field such that $char K\ne p$ and $K$ contains a
primitive $p^e$-th root of unity, and let $\zeta$ be a primitive
$p^m$-th root of unity. Then $K(x_0,x_1,\dots,x_{p^n-1})^G$ is
rational over $K$, where $G$ acts on $x_0,\dots,x_{p^n-1}$ by
\begin{eqnarray*}
\sigma&:&x_i\mapsto \zeta^{k^i}x_i,\\
\tau&:&x_0\mapsto x_1\mapsto\cdots\mapsto x_{p^n-1}\mapsto x_0.
\end{eqnarray*}
\end{t2.6}

\newtheorem{t2.9}[t2.1]{Theorem}
\begin{t2.9}\label{t2.9}
{\rm (}\cite{Ha1,Ha2}{\rm )} Let $K$ be any field and $G$ a finite
subgroup of $\GL(2,\mathbb Z)$. Then the fixed field $K(x, y)^G$
under monomial action of $G$ is rational over $K$.
\end{t2.9}

Finally, we give a Lemma, which can be extracted from some proofs in
\cite{Ka2,HuK}.

\newtheorem{l2.7}[t2.1]{Lemma}
\begin{l2.7}\label{l2.7}
Let $\tau$ be a cyclic group of order $n>1$, acting on
$L(v_1,\dots,v_{n-1})$, the rational function field of $n-1$
variables over a field $L$ such that
\begin{eqnarray*}
\tau&:&v_1\mapsto v_2\mapsto\cdots\mapsto v_{n-1}\mapsto (v_1\cdots
v_{n-1})^{-1}\mapsto v_1.
\end{eqnarray*}
If $L$ contains a primitive $n$th root of unity $\xi$, then
$K(v_1,\dots,v_{n-1})=K(s_1,\dots,s_{n-1})$ where $\tau:s_i\mapsto
\xi^is_i$ for $1\leq i\leq n-1$.
\end{l2.7}
\begin{proof}
Define $w_0=1+v_1+v_1v_2+\cdots+v_1v_2\cdots
v_{n-1},w_1=(1/w_0)-1/n,w_{i+1}=(v_1v_2\cdots v_i/w_0)-1/n$ for
$1\leq i\leq n-1$. Thus $K(v_1,\dots,v_{n-1})=K(w_1,\dots,w_n)$ with
$w_1+w_2+\cdots+w_n=0$ and
\begin{eqnarray*}
\tau&:&w_1\mapsto w_2\mapsto\cdots\mapsto w_{n-1}\mapsto w_n\mapsto
w_1.
\end{eqnarray*}
Define $s_i=\sum_{1\leq j\leq n}\xi^{-ij}w_j$ for $1\leq i\leq n-1$.
Then $K(w_1,\dots,w_n)=K(s_1,\dots,s_{n-1})$ and $\tau:s_i\mapsto
\xi^is_i$ for $1\leq i\leq n-1$.
\end{proof}

\section{The cyclic extensions of the metacyclic $p$-groups}
\label{3}

We are going to consider first the central cyclic extensions of the
abelian metacyclic $p$-groups. Let $G$ be an abelian metacyclic
$p$-group generated by two elements $\sigma$ and $\tau$ with
relations $\sigma^{p^a}=\tau^{p^b}=1$ and
$\tau^{-1}\sigma\tau=\sigma$. Let $\widetilde G$ be a central cyclic
extension of $G$ by $C_{p^t}$ such that the pre-image of
$[\sigma,\tau]=\sigma^{-1}\tau^{-1}\sigma\tau$ in $\widetilde G$ is
of order $p^t$.  Then we have the group extension
\begin{equation}\label{e3.1}
1\longrightarrow C_{p^t}\longrightarrow \widetilde G\longrightarrow
G\cong C_{p^a}\times C_{p^b} \longrightarrow 1,
\end{equation}
where $C_{p^t}=\langle\rho\rangle\leq Z(\widetilde G),t\leq b\leq a$
and the pre-image of $[\sigma,\tau]=\sigma^{-1}\tau^{-1}\sigma\tau$
in $\widetilde G$ is of order $p^t$. If we put $\vert\widetilde
G\vert=p^n$, we get $n=a+b+t$. According to \cite{AMM} any
$2$-generator group of order $p^n$ and nilpotency class $2$ is a
central extension of the form \eqref{e3.1}. Moreover, for any
positive partition $(a,b,t)$ of $n$, the set of nonisomorphic
central extensions of the form \eqref{e3.1} with nilpotency class
exactly $2$ is nonempty \cite[Lemma 2.2]{AMM}.

Any group $\widetilde G$ then has the presentation
\begin{equation*}
\widetilde
G=\langle\sigma,\tau,\rho:\sigma^{p^a}=\rho^{ip^\alpha},\tau^{p^b}=\rho^{jp^\beta},\rho^{p^t}=1,[\sigma,\tau]=\rho^l,
\rho - \text{central}\rangle,
\end{equation*}
where $i,j,l$ are positive integers, $0\leq i,j,l<p^t$,
$\gcd(ijl,p)=1,0\leq\alpha,\beta\leq t;a\geq b\geq t\geq 1$ and
$a+b+t=n$. Note that the commutator subgroup $\widetilde G'$ is
cyclic and is generated by $\rho^l$.

Bacon and Kappe \cite{BK} give a classification of these groups, but
with some omissions that have been corrected recently in \cite{AMM}.
We need not the full classification for our purposes, so we will not
give it in our paper. Instead, we write the following result from
\cite{AMM} which gives a more convenient form of the above
presentations.

\newtheorem{p3.1}{Proposition}[section]
\begin{p3.1}\label{p3.1}
{\rm (}\cite[Proposition 3.1]{AMM}{\rm )} Fix $a\geq b\geq t\geq 1$,
let $0\leq\alpha,\beta\leq t$, and let $i,j,l$ be positive integers,
$0\leq i,j,l<p^t$, with $\gcd(ijl,p)=1$. Then the groups
\begin{equation*}
\widetilde
G=\langle\sigma,\tau,\rho:\sigma^{p^a}=\rho^{ip^\alpha},\tau^{p^b}=\rho^{jp^\beta},\rho^{p^t}=1,[\sigma,\tau]=\rho^l,
\rho - \text{central}\rangle,
\end{equation*}
and
\begin{equation*}
\widetilde
H=\langle\sigma,\tau,\rho:\sigma^{p^a}=\rho^{p^\alpha},\tau^{p^b}=\rho^{p^\beta},\rho^{p^t}=1,[\sigma,\tau]=\rho,
\rho - \text{central}\rangle,
\end{equation*}
are isomorphic.
\end{p3.1}

Now, let $G$ be any nonabelian metacyclic $p$-group generated by two
elements $\sigma$ and $\tau$ with relations
$\sigma^{p^a}=1,\tau^{p^b}=\sigma^{p^c}$ and
$\tau^{-1}\sigma\tau=\sigma^k$ for $k=\varepsilon+p^r$ where
$\varepsilon=1$ if $p$ is odd, $\varepsilon=\pm 1$ if $p=2$, and
$a,b,c,r\geq 0$ are subject to some restrictions. For example, we
have $a,b,c\geq r$. For the description of these restrictions see
e.g. \cite[p. 564]{Ka1}.

Let $t$ be a positive integer such that $t\leq \min\{a,b\}$ and let
$\widetilde G$ be a central cyclic extension of a nonabelian
metacyclic $p$-group $G$ by $C_{p^t}$ such that the pre-image of
$\sigma^{-(k-1)}[\sigma,\tau]=\sigma^{-k}\tau^{-1}\sigma\tau$ in
$\widetilde G$ is of order $p^t$.  Then we have the group extension
\begin{equation}\label{e3.2}
1\longrightarrow C_{p^t}\longrightarrow \widetilde G\longrightarrow
G\longrightarrow 1,
\end{equation}
where $C_{p^t}=\langle\rho\rangle\leq Z(\widetilde G)$.

The group $\widetilde G$ then has the following presentation:
\begin{equation*}
\widetilde
G=\langle\sigma,\tau,\rho:\sigma^{p^a}=\rho^{ip^\alpha},\tau^{p^b}=\sigma^{p^c}\rho^{jp^\beta},\rho^{p^t}=1,\tau^{-1}\sigma\tau=\sigma^k\rho^l,
\rho - \text{central}\rangle,
\end{equation*}
where $i,j,l$ are positive integers, $0\leq i,j,l<p^t$,
$\gcd(ijl,p)=1,0\leq\alpha,\beta\leq t$ and $k=\varepsilon+p^r$.

\newtheorem{p3.2}[p3.1]{Proposition}
\begin{p3.2}\label{p3.2}
Fix $a\geq t, b\geq t, 0\leq c\leq a,r\leq\min\{a,b,c\}$. Let
$0\leq\alpha,\beta\leq t$, and let $i,j,l$ be positive integers,
$0\leq i,j,l<p^t$, with $\gcd(ijl,p)=1$. Then there exist positive
integers $m,n:0<m<p^a,0<n<p^t,\gcd(mn,p)=1$ such that the group
\begin{equation*}
\widetilde
G=\langle\sigma,\tau,\rho:\sigma^{p^a}=\rho^{ip^\alpha},\tau^{p^b}=\sigma^{p^c}\rho^{jp^\beta},\rho^{p^t}=1,\tau^{-1}\sigma\tau=\sigma^k\rho^l,
\rho - \text{central}\rangle,
\end{equation*}
is isomorphic to the group
\begin{equation}\label{e3.3}
\widetilde
H=\langle\sigma,\tau,\rho:\sigma^{p^a}=\rho^{np^\alpha},\tau^{p^b}=\sigma^{mp^c}\rho^{p^\beta},\rho^{p^t}=1,\tau^{-1}\sigma\tau=\sigma^k\rho,
\rho - \text{central}\rangle.
\end{equation}
\end{p3.2}
\begin{proof}
If $\beta=t$, set $m=1,\sigma_1=\sigma,\tau_1=\tau$. Choose an
integer $l_1$ such that $ll_1\equiv 1\pmod{p^t}$, i.e.,
$\rho=(\rho^l)^{l_1}$. Define $n=il_1$ and $\rho_1=\rho^l$. Then the
elements $\sigma_1,\tau_1,\rho_1$ of $\widetilde G$ satisfy the same
relations as $\sigma,\tau,\rho\in\widetilde H$. Whence we have an
onto homomorphism $H\to G$ that maps
$\sigma\mapsto\sigma_1,\tau\mapsto\tau_1$ and $\rho\mapsto\rho_1$.
Since the two groups have the same order, this map is an
isomorphism.

If $\beta<t$, then pick $s$ such that $ls\equiv
j\pmod{p^{t-\beta}}$, and set $\sigma_1=\sigma^s,\tau_1=\tau$, and
$\rho_1=\rho^{ls}$. Then
$\tau_1^{-1}\sigma_1\tau_1=\sigma_1^k\rho_1,\tau_1^{p^b}=\sigma^{p^c}\rho^{jp^\beta}=\sigma_1^{mp^c}\rho^{lsp^\beta}=\sigma_1^{mp^c}\rho_1^{p^\beta}$
for $m$ such that $ms\equiv 1\pmod{\ord(\sigma)}$ (i.e.,
$\sigma=\sigma_1^m$); and $\sigma_1^{p^a}=\rho_1^{np^\alpha}$, where
$n=isu$ for some $u$ such that $uls\equiv 1\pmod{p^t}$ (i.e.,
$\rho=\rho_1^u$). Again, we obtain a homomorphism from $\widetilde
H$ onto $\widetilde G$, showing that $\widetilde G\cong\widetilde
H$.
\end{proof}

We need to calculate the commutator subgroup $\widetilde G'$ in
order to apply Theorem \ref{t2.8} effectively.

\newtheorem{p3.3}[p3.1]{Proposition}
\begin{p3.3}\label{p3.3}
Let the group $\widetilde G$ be isomorphic to the group $\widetilde
H$ with a presentation of the form \eqref{e3.3}. The commutator
subgroup $\widetilde H'$ is cyclic and is generated by
$\sigma^{k-1}\rho$.
\end{p3.3}
\begin{proof}
First, note that $\tau^{-i}\sigma^j\tau^i=\sigma^{jk^i}\rho^{j\cdot
w_i}$ for $1\leq i\leq \ord(\tau),1\leq j\leq\ord(\sigma)$ and
$w_i=1+k+k^2+\cdots+k^{i-1}$.

Let us calculate now an arbitrary commutator:
{\allowdisplaybreaks\begin{eqnarray*}
A&=&[\tau^i\sigma^j\rho^l,\tau^u\sigma^v\rho^w]=[\tau^i\sigma^j,\tau^u\sigma^v]=(\tau^i\sigma^j)^{-1}(\tau^u\sigma^v)^{-1}(\tau^i\sigma^j)(\tau^u\sigma^v)\\
&=&\sigma^{-j}(\tau^{-i}\sigma^{-v}\tau^i)(\tau^{-u}\sigma^j\tau^u)\sigma^v=\sigma^{-j}(\sigma^{-vk^i}\rho^{-vw_i})(\sigma^{jk^u}\rho^{jw_u})\sigma^v\\
&=&\sigma^B\rho^C\quad \text{for}\\
B&=&-j+jk^u-vk^i+v,\\
C&=&-vw_i+jw_u.
\end{eqnarray*}}
We have now $B=j(k^u-1)-v(k^i-1)=j(k-1)w_u-v(k-1)w_i=(k-1)C$, so
$A=(\sigma^{k-1}\rho)^C$. We are done.
\end{proof}

In the proof of Theorem \ref{t1.5} we are going to consider
different cases which appear according to the values of
$\alpha,\beta,c$ and $\varepsilon$. Namely, the group $\widetilde G$
will be equal to one of the following $16$ groups.
\begin{enumerate}
\item $\varepsilon=1,c=a,\alpha=\beta=t$\\
$G_1 =
\langle\sigma,\tau,\rho:\sigma^{p^a}=\tau^{p^b}=\rho^{p^t}=1,\tau^{-1}\sigma\tau=\sigma^{1+p^r}\rho\rangle$,
\item $\varepsilon=1,c=a,\alpha=t,\beta<t$\\
$G_2 =
\langle\sigma,\tau,\rho:\sigma^{p^a}=\rho^{p^t}=1,\tau^{p^b}=\rho^{p^\beta},\tau^{-1}\sigma\tau=\sigma^{1+p^r}\rho\rangle$,
\item $\varepsilon=1,c=a,\alpha<t,\beta=t$\\
$G_3 =
\langle\sigma,\tau,\rho:\sigma^{p^a}=\rho^{p^\alpha},\tau^{p^b}=\rho^{p^t}=1,\tau^{-1}\sigma\tau=\sigma^{1+p^r}\rho\rangle$,
\item $\varepsilon=1,c=a,\alpha<t,\beta<t$\\
$G_4 =
\langle\sigma,\tau,\rho:\sigma^{p^a}=\rho^{p^\alpha},\tau^{p^b}=\rho^{p^\beta},\rho^{p^t}=1,\tau^{-1}\sigma\tau=\sigma^{1+p^r}\rho\rangle$,
\item $\varepsilon=-1,c=a,\alpha=\beta=t$\\
$G_5 =
\langle\sigma,\tau,\rho:\sigma^{2^a}=\tau^{2^b}=\rho^{2^t}=1,\tau^{-1}\sigma\tau=\sigma^{-1+2^r}\rho\rangle$,
\item $\varepsilon=-1,c=a,\alpha=t,\beta<t$\\
$G_6 =
\langle\sigma,\tau,\rho:\sigma^{2^a}=\rho^{2^t}=1,\tau^{2^b}=\rho^{2^\beta},\tau^{-1}\sigma\tau=\sigma^{-1+2^r}\rho\rangle$,
\item $\varepsilon=-1,c=a,\alpha<t,\beta=t$\\
$G_7 =
\langle\sigma,\tau,\rho:\sigma^{2^a}=\rho^{2^\alpha},\tau^{2^b}=\rho^{2^t}=1,\tau^{-1}\sigma\tau=\sigma^{-1+2^r}\rho\rangle$,
\item $\varepsilon=-1,c=a,\alpha<t,\beta<t$\\
$G_8 =
\langle\sigma,\tau,\rho:\sigma^{2^a}=\rho^{2^\alpha},\tau^{2^b}=\rho^{2^\beta},\rho^{2^t}=1,\tau^{-1}\sigma\tau=\sigma^{-1+2^r}\rho\rangle$,
\item $\varepsilon=1,c<a,\alpha=\beta=t$\\
$G_9 =
\langle\sigma,\tau,\rho:\sigma^{p^a}=\rho^{p^t}=1,\tau^{p^b}=\sigma^{p^c},\tau^{-1}\sigma\tau=\sigma^{1+p^r}\rho\rangle$,
\item $\varepsilon=1,c<a,\alpha=t,\beta<t$\\
$G_{10} =
\langle\sigma,\tau,\rho:\sigma^{p^a}=\rho^{p^t}=1,\tau^{p^b}=\sigma^{p^c}\rho^{p^\beta},\tau^{-1}\sigma\tau=\sigma^{1+p^r}\rho\rangle$,
\item $\varepsilon=1,c<a,\alpha<t,\beta=t$\\
$G_{11} =
\langle\sigma,\tau,\rho:\sigma^{p^a}=\rho^{p^\alpha},\tau^{p^b}=\sigma^{p^c},\rho^{p^t}=1,\tau^{-1}\sigma\tau=\sigma^{1+p^r}\rho\rangle$,
\item $\varepsilon=1,c<a,\alpha<t,\beta<t$\\
$G_{12} =
\langle\sigma,\tau,\rho:\sigma^{p^a}=\rho^{p^\alpha},\tau^{p^b}=\sigma^{p^c}\rho^{p^\beta},\rho^{p^t}=1,\tau^{-1}\sigma\tau=\sigma^{1+p^r}\rho\rangle$,
\item $\varepsilon=-1,c<a,\alpha=\beta=t$\\
$G_{13} =
\langle\sigma,\tau,\rho:\sigma^{2^a}=\rho^{2^t}=1,\tau^{2^b}=\sigma^{2^c},\tau^{-1}\sigma\tau=\sigma^{-1+2^r}\rho\rangle$,
\item $\varepsilon=-1,c<a,\alpha=t,\beta<t$\\
$G_{14} =
\langle\sigma,\tau,\rho:\sigma^{2^a}=\rho^{2^t}=1,\tau^{2^b}=\sigma^{2^c}\rho^{2^\beta},\tau^{-1}\sigma\tau=\sigma^{-1+2^r}\rho\rangle$,
\item $\varepsilon=-1,c<a,\alpha<t,\beta=t$\\
$G_{15} =
\langle\sigma,\tau,\rho:\sigma^{2^a}=\rho^{2^\alpha},\tau^{2^b}=\sigma^{2^c},\rho^{2^t}=1,\tau^{-1}\sigma\tau=\sigma^{-1+2^r}\rho\rangle$,
\item $\varepsilon=-1,c<a,\alpha<t,\beta<t$\\
$G_{16} =
\langle\sigma,\tau,\rho:\sigma^{2^a}=\rho^{2^\alpha},\tau^{2^b}=\sigma^{2^c}\rho^{2^\beta},\rho^{2^t}=1,\tau^{-1}\sigma\tau=\sigma^{-1+2^r}\rho\rangle$,
\end{enumerate}

\section{Proof of Theorem \ref{t1.4}}
\label{4}

The case (i) follows from Theorem \ref{t2.3}.

(ii) We will divide the proof into two steps.

{\bf Step I.} Let $C=C_{p^{\alpha_1}}\times\cdots\times
C_{p^{\alpha_s}}\leq Z(\widetilde G)$. Denote by $\sigma$ and $\tau$
the preimages of the generators of $G$ and by $\rho_1,\dots,\rho_s$
the generators of $C$, i.e., $\rho_i^{p^{\alpha_i}}=1$. Then
$[\sigma,\tau]=\rho_1^{\beta_1}\cdots\rho_s^{\beta_s}$ for
$\beta_i\geq 0$, and we can assume for abuse of notation that
$p^t=\ord(\rho_1^{\beta_1})=\max\{\ord(\rho_1^{\beta_1}),\dots,\ord(\rho_s^{\beta_s})\}$.
It follows that $\widetilde
G'\cap\langle\rho_2,\dots,\rho_s\rangle=\{1\}$. Now we can apply
Theorem \ref{t2.8} reducing the rationality problem of $K(\widetilde
G)$ over $K$ to the rationality problem of $K(\widetilde
G/\langle\rho_2,\dots,\rho_k\rangle)$ over $K$. Note that
$\exp(\langle\rho_2,\dots,\rho_k\rangle)\leq p^t\leq
p^a=\exp(\widetilde G_1/\widetilde G_1')$ for $\widetilde
G_1=\widetilde G/\langle\rho_2,\dots,\rho_k\rangle$.

Therefore, we will assume that $\widetilde G$ has the presentation
\begin{equation*}
\widetilde
G=\langle\sigma,\tau,\rho:\sigma^{p^a}=\rho^{ip^\alpha},\tau^{p^b}=\rho^{jp^\beta},\rho^{p^t}=1,[\sigma,\tau]=\rho^l,
\rho - \text{central}\rangle,
\end{equation*}
where $i,j,l$ are positive integers, $0\leq i,j,l<p^t$,
$\gcd(ijl,p)=1,0\leq\alpha,\beta\leq t;a\geq b\geq t\geq 1$. Note
that the commutator subgroup $\widetilde G'$ is cyclic and is
generated by $\rho^l$.

Moreover, from Proposition \ref{p3.1} it follows that we can assume
that $\widetilde G$ is $2$-generator $p$-group of nilpotency class
$2$ given by the presenation
\begin{equation*}
\widetilde
G=\langle\sigma,\tau,\rho:\sigma^{p^a}=\rho^{p^\alpha},\tau^{p^b}=\rho^{p^\beta},\rho^{p^t}=1,[\sigma,\tau]=\rho,
\rho - \text{central}\rangle,
\end{equation*}
where $0\leq\alpha,\beta\leq t$ and $a\geq b\geq t$.

Define $m=a+t-\alpha$ and $n=b+t-\beta$. Then $p^m=\ord(\sigma)$ and
$p^n=\ord(\tau)$. According to \cite{AMM} the center of $\widetilde
G$ is
\begin{equation*}
Z(\widetilde G)=\{\sigma^{ip^t}\tau^{jp^t}\rho^k: 1\leq i\leq
p^{m-t},1\leq j\leq p^{n-t},1\leq k\leq p^t\}.
\end{equation*}
Consider the following two possibilities for $\alpha$.

{\bf Case 1.} Let $\alpha=t$, i.e., $\sigma^{p^a}=1$. Then
$\widetilde G'\cap \langle\sigma^{p^t}\rangle=\{1\}$, so from
Theorem \ref{t2.8} it follows that we can reduce the rationality
problem to one related to the group generated by elements
$\sigma,\tau$ and $\rho$ such that
$\sigma^{p^t}=1,\tau^{p^b}=\rho^{p^\beta},\rho^{p^t}=1,[\sigma,\tau]=\rho,
\rho - \text{central}$.

{\bf Case 2.} Let $\alpha<t$, i.e., $\sigma^{p^a}=\rho^{p^\alpha}\ne
1$. If $a=t$, we obtain that $\widetilde G$ is generated by elements
$\sigma,\tau$ and $\rho$ such that
$\sigma^{p^t}=\rho^{p^\alpha},\tau^{p^b}=\rho^{p^\beta},\rho^{p^t}=1,[\sigma,\tau]=\rho,
\rho - \text{central}$. Let $a>t$. (Recall that $a\geq b\geq t$.)
Then $\sigma^{p^{a-1}}=\sigma^{p^{a-t-1}\cdot p^t}\in Z(\widetilde
G)$ and we get $(\sigma^{p^{a-1}}\rho^{-p^{\alpha-1}})^p=1$. Then
$\widetilde G'\cap
\langle\sigma^{p^{a-1}}\rho^{-p^{\alpha-1}}\rangle=\{1\}$, so we can
apply Theorem \ref{t2.8} reducing the rationality problem of
$K(\widetilde G)$ to the rationality problem of $K(\widetilde G_1)$
over $K$, where $\widetilde G_1\cong\widetilde
G/\langle\sigma^{p^{a-1}}\rho^{-p^{\alpha-1}}\rangle$. The group
$\widetilde G_1$ is generated by elements $\sigma,\tau$ and $\rho$
such that
$\sigma^{p^{a-1}}=\rho^{p^{\alpha-1}},\tau^{p^b}=\rho^{p^\beta},\rho^{p^t}=1,[\sigma,\tau]=\rho,
\rho - \text{central}$. We can proceed in the same way until we
obtain either a group $\widetilde G_k$ such that
$\sigma^{p^t}=\rho^{p^{\alpha-u}}$ for $u=a-t$ and $\alpha\geq u$;
or a metacyclic $p$-group such that $\sigma^{p^{a-l}}=\rho$ for some
$l$. For the second group we can apply Theorem \ref{t1.2}.

Similarly, we can consider the cases $\beta=t$ and $\beta<t$. In
this way, we may assume that $\widetilde G$ is generated by elements
$\sigma,\tau$ and $\rho$ such that
$\sigma^{p^t}=\rho^{p^\alpha},\tau^{p^t}=\rho^{p^\beta},\rho^{p^t}=1,[\sigma,\tau]=\rho,
\rho - \text{central}$. Note that anytime we applied Theorem
\ref{t2.8} so far, the condition for the root of unity is satisfied,
since $\exp(\widetilde G/\widetilde G')\leq p^a$.

{\bf Step II.} We are going to consider the regular representation
of $\widetilde G$. Let $V$ be a $K$-vector space whose dual space
$V^*$ is defined as $V^*=\bigoplus_{g\in\widetilde G}K\cdot x(g)$
where $\widetilde G$ acts on $V^*$ by $h\cdot x(g)=x(hg)$ for any
$h,g\in\widetilde G$. Thus $K(V)^{\widetilde
G}=K(x(g):g\in\widetilde G)^{\widetilde G}=K(\widetilde G)$.

For $m=2t-\alpha$ and $n=2t-\beta$ define $\zeta=\zeta_{p^m}$, a
primitive $p^m$-th root of unity; $\eta=\zeta^{p^{m-n}}$, a
primitive $p^n$-th root of unity; and $\xi=\zeta^{p^{t-\alpha}}$, a
primitive $p^t$-th root of unity. Since $K$ contains a primitive
$p^{2t}$-th root of unity, we have $\zeta,\eta,\xi\in K$. For $0\leq
i\leq p^t-1$ define $x_i\in V^*$ by
\begin{equation*}
x_i=\sum_{j,k}\eta^{-j-kp^{b-\beta}}x(\sigma^i\tau^j\rho^k),
\end{equation*}
where $0\leq j\leq p^t-1,0\leq k\leq p^t-1$. The actions of $\sigma,
\tau$ and $\rho$ are given by {\allowdisplaybreaks
\begin{align*}
\sigma\ :\ &x_0\mapsto x_1\mapsto\cdots\mapsto x_{p^t-1}\mapsto\xi^{p^\alpha} x_0,\\
\tau\ :\ &x_i\mapsto \eta\xi^ix_i,\\
\rho\ :\ &x_i\mapsto \xi x_i,
\end{align*}}
for $0\leq i\leq p^t-1$. We find that $Y=\bigoplus_{0\leq i\leq
p^t-1}K\cdot x_i$ is a faithful $\widetilde G$-subspace of $V^*$.
Thus, by Theorem \ref{t2.1}, it suffices to show that $K(x_i:0\leq
i\leq p^t-1)^{\widetilde G}$ is rational over $K$.

Define
\begin{equation*}
y_0=x_0^{p^t},y_1=x_1/x_0,y_2=x_2/x_1,\dots,y_{p^t-1}=x_{p^t-1}/x_{p^t-2}.
\end{equation*}
We have now
$K(x_0,\dots,x_{p^t-1})^{\langle\rho\rangle}=K(y_0,\dots,y_{p^t-1})$
and {\allowdisplaybreaks\begin{eqnarray*}
\sigma&:&y_0\mapsto y_1^{p^t}y_0,y_1\mapsto y_2\mapsto\cdots\mapsto y_{p^t-1}\mapsto \xi^{p^\alpha}/(y_1\cdots y_{p^t-1})\mapsto y_1,\\
\tau&:&y_0\mapsto \eta^{p^t}y_0, y_i\mapsto\xi y_i,
\end{eqnarray*}}
for $1\leq i\leq p^t-1$. From Theorem \ref{t2.2} follows that if
$K(y_1,\dots,y_{p^t-1})^{\langle\sigma,\tau\rangle}$ is rational
over $K$, so is
$K(y_0,\dots,y_{p^t-1})^{\langle\sigma,\tau\rangle}$.

Define
\begin{equation*}
z_2=y_2/y_1,z_3=y_3/y_2,\dots,z_{p^t-1}=y_{p^t-1}/y_{p^t-2}.
\end{equation*}
We have {\allowdisplaybreaks\begin{eqnarray*} \sigma&:&y_1\mapsto
z_2y_1,z_2\mapsto z_3\mapsto\cdots\mapsto z_{p^t-1}\mapsto
\xi^{p^\alpha}/(y_1y_2\cdots y_{p^t-1}^2)=\xi^{p^\alpha}/(y_1^{p^t}z_2^{p^t-1}\cdots z_{p^t-1}^2),\\
\tau&:&y_1\mapsto \xi y_1, z_i\mapsto z_i,
\end{eqnarray*}}
for $2\leq i\leq p^t-1$.

Define $z_1=y_1^{p^t}\xi^{-p^\alpha}$, i.e.,
$y_1^{p^t}=z_1\xi^{p^\alpha}$. We have now
$K(y_1,\dots,y_{p^t-1})^{\langle\tau\rangle}=K(z_1,\dots,z_{p^t-1})$
and {\allowdisplaybreaks\begin{eqnarray*} \sigma&:&z_1\mapsto
z_2^{p^t}z_1,z_2\mapsto z_3\mapsto\cdots\mapsto z_{p^t-1}\mapsto
1/(z_1z_2^{p^t-1}z_3^{p^t-2}\cdots z_{p^t-1}^2).
\end{eqnarray*}}
Define $s_1=z_2,s_i=\tau^{i-1}\cdot z_2$ for $2\leq i\leq p^t-1$.
Then $K(z_i:1\leq i\leq p^t-1)=K(s_i:1\leq i\leq p^t-1)$ and
\begin{align*}
\tau\ :\ &s_1\mapsto s_2\mapsto\cdots\mapsto s_{p^t-1}\mapsto
(s_1s_2\cdots s_{p^t-1})^{-1}.
\end{align*}
The action of $\tau$ can be linearized according to Lemma
\ref{l2.7}. Thus $K(s_i:1\leq i\leq p^t-1)^{\langle\tau\rangle}$ is
rational over $K$ by Theorem \ref{t1.1}. We are done.

\section{Proof of Theorem \ref{t1.5}}
\label{5}

The case (i) follows from Theorem \ref{t2.3}.

(ii) We will divide the proof into two steps.

{\bf Step I.} Assume that $\widetilde G$ has the following
presentation:
\begin{equation*}
\widetilde
G=\langle\sigma,\tau,\rho:\sigma^{p^a}=\rho^{p^\alpha},\tau^{p^b}=\sigma^{p^c}\rho^{p^\beta},\rho^{p^t}=1,\tau^{-1}\sigma\tau=\sigma^k\rho,
\rho - \text{central}\rangle,
\end{equation*}
where $a\geq t, b\geq t, 0\leq\alpha,\beta\leq t$ and
$k=\varepsilon+p^r$. Note that we have the relations
$\tau^{-i}\sigma\tau^i=\sigma^{k^i}\rho^{w_i}$, where
$w_i=1+k+k^2+\cdots+k^{i-1}$ and $1\leq i\leq\ord(\tau)$.  Clearly,
$\widetilde G$ is isomorphic to some group $\widetilde G_i$ for
$1\leq i\leq 16$ given in Section \ref{3}. We are going to consider
each case separately.

{\bf Case 1.} $\widetilde G=G_1$, where $G_1$ is the group in
Section \ref{3}. From Proposition \ref{p3.3} it follows that
$\widetilde G'$ is cyclic and is generated by $\sigma^{p^r}\rho$.

{\bf Subcase 1.a.} Let $a\geq r+t$. Then
$(\sigma^{p^r}\rho)^{p^{a-r}}=\rho^{p^{a-r}}=1$, so $\widetilde
G'\cap\langle\rho\rangle=1$. Theorem \ref{t2.8} then implies that we
can reduce the rationality problem of $K(\widetilde G)$ to
$K(\widetilde G/\langle\rho\rangle)$ over $K$, where $\widetilde
G/\langle\rho\rangle$ clearly is a metacyclic $p$-group.

{\bf Subcase 1.b.} Let $a<r+t$. Let $V$ be a $K$-vector space whose
dual space $V^*$ is defined as $V^*=\bigoplus_{g\in\widetilde
G}K\cdot x(g)$ where $\widetilde G$ acts on $V^*$ by $h\cdot
x(g)=x(hg)$ for any $h,g\in\widetilde G$. Thus $K(V)^{\widetilde
G}=K(x(g):g\in\widetilde G)^{\widetilde G}=K(\widetilde G)$.

Define $X_1,X_2\in V^*$ by
\begin{equation*}
X_1=\sum_{i=0}^{p^a-1}x(\sigma^i),~ X_2=\sum_{i=0}^{p^t-1}x(\rho^i).
\end{equation*}
Note that $\sigma\cdot X_1=X_1$ and $\rho\cdot X_2=X_2$.

Let $\zeta=\zeta_{p^a}\in K$ be a primitive $p^a$-th root of unity
and define $\xi=\zeta^{p^{a-t}}$. Thus $\xi$ is a primitive $p^t$-th
root of unity. Define $Y_1,Y_2\in V^*$ by
\begin{equation*}
Y_1=\sum_{i=0}^{p^t-1}\xi^{-i}\rho^i\cdot X_1,~
Y_2=\sum_{i=0}^{p^a-1}\zeta^{-i}\sigma^i\cdot X_2.
\end{equation*}

It follows that {\allowdisplaybreaks \begin{align*}
\sigma\ :\ &Y_1\mapsto Y_1,Y_2\mapsto\zeta Y_2,\\
\rho\ :\ &Y_1\mapsto\xi Y_1,Y_2\mapsto Y_2.
\end{align*}}
Thus $K\cdot Y_1+K\cdot Y_2$ is a representation space of the
subgroup $\langle\sigma,\rho\rangle$.

Define $x_i=\tau^i\cdot Y_1,y_i=\tau^i\cdot Y_2$ for $0\leq i\leq
p^b-1$. We have now {\allowdisplaybreaks \begin{align*}
\sigma\ :\ &x_i\mapsto \xi^{w_i} x_i,~ y_i\mapsto\zeta^{k^i}y_i\\
\tau\ :\ &x_0\mapsto x_1\mapsto\cdots\mapsto x_{p^b-1}\mapsto x_0,\\
&y_0\mapsto y_1\mapsto\cdots\mapsto y_{p^b-1}\mapsto y_0,\\
\rho\ :\ &x_i\mapsto\xi x_i,~ y_i\mapsto y_i,
\end{align*}}
for $0\leq i\leq p^b-1$.

We find that $Y=(\bigoplus_{0\leq i\leq p^b-1}K\cdot
x_i)\oplus(\bigoplus_{0\leq i\leq p^b-1}K\cdot y_i)$ is a faithful
$\widetilde G$-subspace of $V^*$. Thus, by Theorem \ref{t2.1}, it
suffices to show that $K(x_i,y_i:0\leq i\leq p^b-1)^{\widetilde G}$
is rational over $K$.

For $1\leq i\leq p^b-1$, define $u_i=x_i/x_{i-1}$ and
$v_i=y_i/y_{i-1}$. Thus $K(x_i,y_i:0\leq i\leq
p^b-1)=K(x_0,y_0,u_i,v_i:1\leq i\leq p^b-1)$ and for every
$g\in\widetilde G$
\begin{equation*}
g\cdot x_0\in K(u_i,v_i:1\leq i\leq p^t-1)\cdot x_0,~ g\cdot y_0\in
K(u_i,v_i:1\leq i\leq p^b-1)\cdot y_0,
\end{equation*}
while the subfield $K(u_i,v_i:1\leq i\leq p^b-1)$ is invariant by
the action of $\widetilde G$. Thus $K(x_i,y_i:0\leq i\leq
p^b-1)^{\widetilde G}=K(u_i,v_i:1\leq i\leq p^b-1)^{\widetilde
G}(u,v)$ for some $u,v$ such that $\sigma(v)=\tau(v)=\rho(v)=v$ and
$\sigma(u)=\tau(u)=\rho(u)=u$. We have now {\allowdisplaybreaks
\begin{align*}
\sigma\ :\ &u_i\mapsto \xi^{k^{i-1}} u_i,~ v_i\mapsto\zeta^{k^i-k^{i-1}}v_i,~ v\mapsto v\\
\tau\ :\ &u_1\mapsto u_2\mapsto\cdots\mapsto u_{p^b-1}\mapsto (u_1u_2\cdots u_{p^b-1})^{-1},\\
&v_1\mapsto v_2\mapsto\cdots\mapsto v_{p^b-1}\mapsto (v_1v_2\cdots v_{p^b-1})^{-1},\\
&~ v\mapsto v,\\
\rho\ :\ &u_i\mapsto u_i,~ v_i\mapsto v_i,~ v\mapsto v
\end{align*}}
for $1\leq i\leq p^b-1$. From Theorem \ref{t2.2} it follows that if
$K(u_i,v_i:1\leq i\leq p^b-1)^{\widetilde G}(v)$ is rational over
$K$, so is $K(x_i,y_i:0\leq i\leq p^b-1)^{\widetilde G}$ over $K$.

Since $\rho$ acts trivially on $K(u_i,v_i:1\leq i\leq p^t-1)$, we
find that $K(u_i,v_i:1\leq i\leq p^t-1)^{\widetilde
G}=K(u_i,v_i:1\leq i\leq p^t-1)^{\langle\sigma,\tau\rangle}$.

Recall that $a<r+t$. So we can write $r=a-t+r_1$ for some $r_1\geq
1$. It follows that
$\zeta^{k^{i-1}(k-1)}=\zeta^{k^{i-1}p^{a-t+r_1}}=\xi^{k^{i-1}p^{r_1}}$
for $1\leq i\leq p^b-1$.

Define $w_i=v_i/(u_i)^{p^{r_1}}$. We have now that $K(u_i,v_i:1\leq
i\leq p^b-1)=K(u_i,w_i:1\leq i\leq p^b-1)$ and {\allowdisplaybreaks
\begin{align*}
\sigma\ :\ &u_i\mapsto \xi^{k^{i-1}} u_i,~ w_i\mapsto w_i,~ v\mapsto v\\
\tag{5.1} \tau\ :\ &u_1\mapsto u_2\mapsto\cdots\mapsto u_{p^b-1}\mapsto (u_1u_2\cdots u_{p^b-1})^{-1},\\
&w_1\mapsto w_2\mapsto\cdots\mapsto w_{p^b-1}\mapsto (w_1w_2\cdots
w_{p^b-1})^{-1},\\
&v\mapsto v.
\end{align*}}
According to Lemma \ref{l2.7}, we can linearize the action of $\tau$
on $w_1,\dots,w_{p^b-1}$.

Write $L=K(v,u_i:1\leq i\leq p^b-1)$ and consider $L(w_i:1\leq i\leq
p^b-1)^{\langle\sigma,\tau\rangle}$. Note that the group
$\langle\sigma,\tau\rangle$ acts on the field $L(w_i)$ as
$\langle\sigma,\tau\rangle/\langle\sigma^{p^t}\rangle$ and is
faithful on $L$. Thus we may apply Theorem \ref{t2.1} to
$L(w_i:1\leq i\leq p^b-1)^{\langle\sigma,\tau\rangle}$. It remains
to show that $L^{\langle\sigma,\tau\rangle}$ is rational over $K$.

Let $\eta$ be a primitive $p^{r+t}$-th root of unity such that
$\xi=\eta^{p^r}$. Whence $\xi^{k^{i-1}}=\eta^{k^i-k^{i-1}}$.

Consider the metacyclic $p$-group $\widetilde
G_1=\langle\sigma,\tau:\sigma^{p^{r+t}}=\tau^{p^b}=1,\tau^{-1}\sigma\tau=\sigma^k\rangle$.

Define $X=\sum_{0\leq j\leq
p^{r+t}-1}\eta^{-j}x(\sigma^j),V_i=\tau^i X$ for $0\leq i\leq
p^t-1$. It follows that
\begin{eqnarray*}
\sigma&:&V_i\mapsto \eta^{k^i}V_i,\\
\tau&:&V_0\mapsto V_1\mapsto\cdots\mapsto V_{p^b-1}\mapsto V_0.
\end{eqnarray*}
Note that $K(V_0,V_1,\dots,V_{p^b-1})^{\widetilde G_1}$ is rational
by Theorem \ref{t2.6}.

Define $v_i=V_i/V_{i-1}$ for $1\leq i\leq p^b-1$. Then
$K(V_0,V_1,\dots,V_{p^b-1})^{\widetilde G_1}=K(v_1,v_2,\dots,$
$v_{p^b-1})^{\widetilde G_1}(V)$ where
\begin{eqnarray*}
\sigma&:&v_i\mapsto \eta^{k^i-k^{i-1}}v_i,~ V\mapsto V\\
\tau&:&v_1\mapsto v_2\mapsto\cdots\mapsto v_{p^b-1}\mapsto
(v_1v_2\cdots v_{p^b-1})^{-1},~ V\mapsto V.
\end{eqnarray*}

Whence $K(u_1,\dots,u_{p^b-1})^{\widetilde G}(v)\cong
K(v_1,\dots,v_{p^b-1})^{\widetilde
G_1}(V)=K(V_0,V_1,\dots,V_{p^b-1})^{\widetilde G_1}$ is rational
over $K$.

{\bf Case 2.} $\widetilde G=G_2$. The proof is almost the same as in
Case 1. Here only the action of $\tau$ on $x_i$'s and $u_i$'s is
changed:
\begin{align*}
\tau\ :\ &x_0\mapsto x_1\mapsto\cdots\mapsto
x_{p^b-1}\mapsto\xi^{p^\beta} x_0.
\end{align*}
and, respectively,
\begin{align*}
\tau\ :\ &u_1\mapsto u_2\mapsto\cdots\mapsto
u_{p^b-1}\mapsto\xi^{p^\beta} (u_1u_2\cdots u_{p^b-1})^{-1}.
\end{align*}

Let $\omega\in K$ be a primitive $p^{b+t-\beta}$-th root of unity
such that $\omega^{p^{b-\beta}}=\xi$. Then
$\xi^{p^\beta}=\omega^{p^b}$. Define $U_i=u_i/\omega$ for $1\leq
i\leq p^b-1$. Then we have
\begin{align*}
\tau\ :\ &U_1\mapsto U_2\mapsto\cdots\mapsto U_{p^b-1}\mapsto
(U_1U_2\cdots U_{p^b-1})^{-1}.
\end{align*}
Apply the proof of Case 1.

{\bf Case 3.} $\widetilde G=G_3$. The subgroup
$H=\langle\sigma,\rho\rangle$ is abelian and has an order $p^{a+t}$.
Put $\rho_1=\sigma^{p^{a-\alpha}}\rho^{-1}$. Then
$H\cong\langle\sigma\rangle\times\langle\rho_1\rangle$, where
$\sigma^{p^{a+t-\alpha}}=\rho_1^{p^\alpha}=1$. We have also
$\rho=\sigma^{p^{a-\alpha}}\rho_1^{-1}$.

Define $X_1,X_2\in V^*$ by
\begin{equation*}
X_1=\sum_ix(\rho_1^i),~ X_2=\sum_ix(\sigma^i).
\end{equation*}
Note that $\sigma\cdot X_2=X_2,\rho_1\cdot X_1=X_1$.

Let $\zeta_1\in K$ be a primitive $p^{a+t-\alpha}$-th root of unity.
Put $\xi=\zeta_1^{p^{a-\alpha}}$, a primitive $p^t$-th root of
unity. Let $\zeta_2$ be any primitive $p^\alpha$-th root of unity.
(Note that we will specify $\zeta_2$ a bit later.)

Define $Y_1,Y_2,Y_3\in V^*$ by
\begin{equation*}
Y_1=\sum_{i=0}^{p^{a+t-\alpha}-1}\zeta_1^{-i}\sigma^i\cdot X_1,~
Y_2=\sum_{i=0}^{p^\alpha-1}\zeta_2^{-i}\rho_1^i\cdot X_3.
\end{equation*}

It follows that {\allowdisplaybreaks \begin{align*}
\sigma\ :\ &Y_1\mapsto\zeta_1 Y_1,~ Y_2\mapsto Y_2,\\
\rho_1\ :\ &Y_1\mapsto Y_1,~ Y_2\mapsto\zeta_2 Y_2,\\
\rho\ :\ &Y_1\mapsto\xi Y_1,~ Y_2\mapsto\zeta_2^{-1} Y_2.
\end{align*}}
Thus $K\cdot Y_1+K\cdot Y_2$ is a representation space of the
subgroup $H$.

Define $x_i=\tau^i\cdot Y_1,y_i=\tau^i\cdot Y_2$ for $0\leq i\leq
p^b-1$. We have now {\allowdisplaybreaks
\begin{align*}
\sigma\ :\ &x_i\mapsto \zeta_1^{k^i}\xi^{w_i} x_i,~ y_i\mapsto\zeta_2^{-w_i}y_i,\\
\tau\ :\ &x_0\mapsto x_1\mapsto\cdots\mapsto x_{p^b-1}\mapsto x_0,\\
&y_0\mapsto y_1\mapsto\cdots\mapsto y_{p^b-1}\mapsto y_0,\\
\rho\ :\ &x_i\mapsto\xi x_i,~ y_i\mapsto\zeta_2^{-1} y_i.
\end{align*}}
for $0\leq i\leq p^b-1$. We find that $Y=(\bigoplus_{0\leq i\leq
p^b-1}K\cdot x_i)\oplus(\bigoplus_{0\leq i\leq p^b-1}K\cdot y_i)$ is
a faithful $\widetilde G$-subspace of $V^*$. Thus, by Theorem
\ref{t2.1}, it suffices to show that $K(x_i,y_i:0\leq i\leq
p^b-1)^{\widetilde G}$ is rational over $K$.

For $1\leq i\leq p^b-1$, define $U_i=x_i/x_{i-1}$ and
$V_i=y_i/y_{i-1}$. Thus $K(x_i,y_i:0\leq i\leq
p^b-1)=K(x_0,y_0,U_i,V_i:1\leq i\leq p^b-1)$ and for every
$g\in\widetilde G$
\begin{equation*}
g\cdot x_0\in K(U_i,V_i)\cdot x_0,~ g\cdot y_0\in K(U_i,V_i)\cdot
y_0,
\end{equation*}
while the subfield $K(U_i,V_i:1\leq i\leq p^b-1)$ is invariant by
the action of $\widetilde G$, i.e., {\allowdisplaybreaks
\begin{align*}
\sigma\ :\ &U_i\mapsto \zeta_1^{k^i-k^{i-1}}\xi^{k^{i-1}} U_i,~ V_i\mapsto\zeta_2^{-k^{i-1}} V_i,\\
\tau\ :\ &U_1\mapsto U_2\mapsto\cdots\mapsto U_{p^b-1}\mapsto (U_1\cdots U_{p^b-1})^{-1},\\
&V_1\mapsto V_2\mapsto\cdots\mapsto V_{p^b-1}\mapsto (V_1\cdots V_{p^b-1})^{-1},\\
\rho\ :\ &U_i\mapsto U_i,~ V_i\mapsto V_i.
\end{align*}}
for $1\leq i\leq p^b-1$. From Theorem \ref{t2.2} it follows that if
$K(U_i,V_i:1\leq i\leq p^b-1)^{\widetilde G}$ is rational over $K$,
so is $K(x_i,y_i:0\leq i\leq p^b-1)^{\widetilde G}$ over $K$.

Since $\rho$ acts trivially on $K(U_i,V_i)$, we find that
$K(U_i,V_i)^{\widetilde G}=K(U_i,V_i)^{\langle\sigma,\tau\rangle}$.

{\bf Subcase 3.a.} Let $a-\alpha\leq r$. Thus we can write
$r=a-\alpha+r_1$ for some $r_1\geq 0$. Define
$\zeta_2=\xi^{(1+p^{r_1})p^{t-\alpha}}$, a primitive $p^\alpha$-th
root of unity. Therefore,
$\zeta_1^{k-1}=\zeta_1^{p^{a-\alpha+r_1}}=\xi^{p^{r_1}}$ and also
$\zeta_1^{k^{i-1}(k-1)}=\xi^{k^{i-1}p^{r_1}}$ for all $i$.

Define $v_i=U_i^{p^{t-\alpha}}V_i$. Since
$\zeta_2=\xi^{(1+p^{r_1})p^{t-\alpha}}$, we have
{\allowdisplaybreaks
\begin{align*}
\sigma\ :\ &U_i\mapsto \xi^{k^{i-1}(1+p^{r_1})} U_i,~ v_i\mapsto v_i,\\
\tag{5.2} \tau\ :\ &U_1\mapsto U_2\mapsto\cdots\mapsto U_{p^b-1}\mapsto (U_1\cdots U_{p^b-1})^{-1},\\
&v_1\mapsto v_2\mapsto\cdots\mapsto v_{p^b-1}\mapsto (v_1\cdots
v_{p^b-1})^{-1}
\end{align*}}
for $1\leq i\leq p^b-1$. Compare Formula (5.2) with Formula (5.1).
They look almost the same. Apply the proof of Case 1.

{\bf Subcase 3.b.} Let $a-\alpha>r$. Let $\eta$ be a primitive
$p^{r+t}$-th root of unity such that $\xi=\eta^{p^r}$. Whence
$\xi^{k^{i-1}}=\eta^{k^i-k^{i-1}}$. Since $a+t-\alpha>r+t$, we get
that $\zeta_1\eta$ is a primitive $p^{a+t-\alpha}$-th root of unity.
Put $\mu=(\zeta_1\eta)^{p^{a+t-\alpha-(\alpha+r)}}$, a primitive
$p^{\alpha+r}$-th root of unity, where $\alpha+r<a\leq a+t-\alpha$.
Now, put $\zeta_2=\mu^{p^r}$, a primitive $p^\alpha$-th root of
unity. Thus $\zeta_2^{k^{i-1}}=\mu^{k^i-k^{i-1}}$ for any $i$.

Define $v_i=U_i^{p^{a+t-\alpha-(\alpha+r)}}V_i$ for $1\leq i\leq
p^b-1$. It follows that {\allowdisplaybreaks
\begin{align*}
\sigma\ :\ &U_i\mapsto (\zeta_1\eta)^{k^i-k^{i-1}} U_i,~ v_i\mapsto v_i,\\
\tag{5.3} \tau\ :\ &U_1\mapsto U_2\mapsto\cdots\mapsto U_{p^b-1}\mapsto (U_1\cdots U_{p^b-1})^{-1},\\
&v_1\mapsto v_2\mapsto\cdots\mapsto v_{p^b-1}\mapsto (v_1\cdots
v_{p^b-1})^{-1}
\end{align*}}
for $1\leq i\leq p^b-1$. Compare Formula (5.3) with Formula (5.1).
Considering that $\xi^{k^{i-1}}=\eta^{k^i-k^{i-1}}$, they look
almost the same. Apply the proof of Case 1.

{\bf Case 4.} $\widetilde G=G_4$. The proof is almost the same as in
Case 3. Here only the action of $\tau$ is changed. As we already saw
in Case 2, this issue is solved with a proper adjustment involving
certain root of unity which is in $K$.

{\bf Case 5.} $\widetilde G=G_5$. From Proposition \ref{p3.3} it
follows that $\widetilde G$ is cyclic and is generated by
$\sigma^{-2+2^r}\rho$. If $a\geq t+1$ then
$(\sigma^{-2+2^r}\rho)^{2^{a-1}}=1$ and whence $\widetilde
G\cap\langle\rho\rangle=\{1\}$. Theorem \ref{t2.8} then implies that
we can reduce the rationality problem of $K(\widetilde G)$ to
$K(\widetilde G/\langle\rho\rangle)$ over $K$, where $\widetilde
G/\langle\rho\rangle$ clearly is a metacyclic $p$-group. Therefore,
we may assume that $a=t$.

Define $X_1,X_2\in V^*$ by
\begin{equation*}
X_1=\sum_{i=0}^{2^t-1}x(\sigma^i),~ X_2=\sum_{i=0}^{2^t-1}x(\rho^i).
\end{equation*}
Note that $\sigma\cdot X_1=X_1$ and $\rho\cdot X_2=X_2$.

Let $\xi$ be a primitive $2^t$-th root of unity. Define $Y_1,Y_2\in
V^*$ by
\begin{equation*}
Y_1=\sum_{i=0}^{2^t-1}\xi^{-i}\rho^i\cdot X_1,~
Y_2=\sum_{i=0}^{2^t-1}\xi^{-i}\sigma^i\cdot X_2.
\end{equation*}

It follows that {\allowdisplaybreaks \begin{align*}
\sigma\ :\ &Y_1\mapsto Y_1,Y_2\mapsto\xi Y_2,\\
\rho\ :\ &Y_1\mapsto\xi Y_1,Y_2\mapsto Y_2.
\end{align*}}
Thus $K\cdot Y_1+K\cdot Y_2$ is a representation space of the
subgroup $\langle\sigma,\rho\rangle$.

Define $x_i=\tau^i\cdot Y_1,y_i=\tau^i\cdot Y_2$ for $0\leq i\leq
2^b-1$. We have now {\allowdisplaybreaks \begin{align*}
\sigma\ :\ &x_i\mapsto \xi^{w_i} x_i,~ y_i\mapsto\xi^{k^i}y_i\\
\tau\ :\ &x_0\mapsto x_1\mapsto\cdots\mapsto x_{2^b-1}\mapsto x_0,\\
&y_0\mapsto y_1\mapsto\cdots\mapsto y_{2^b-1}\mapsto y_0,\\
\rho\ :\ &x_i\mapsto\xi x_i,~ y_i\mapsto y_i,
\end{align*}}
for $0\leq i\leq 2^b-1$.

We find that $Y=(\bigoplus_{0\leq i\leq 2^b-1}K\cdot
x_i)\oplus(\bigoplus_{0\leq i\leq 2^b-1}K\cdot y_i)$ is a faithful
$\widetilde G$-subspace of $V^*$. Thus, by Theorem \ref{t2.1}, it
suffices to show that $K(x_i,y_i:0\leq i\leq 2^b-1)^{\widetilde G}$
is rational over $K$.

For $1\leq i\leq 2^b-1$, define $U_i=x_i/x_{i-1}$ and
$V_i=y_i/y_{i-1}$. Thus $K(x_i,y_i:0\leq i\leq
2^b-1)=K(x_0,y_0,U_i,V_i:1\leq i\leq 2^b-1)$ and for every
$g\in\widetilde G$
\begin{equation*}
g\cdot x_0\in K(U_i,V_i:1\leq i\leq p^t-1)\cdot x_0,~ g\cdot y_0\in
K(U_i,V_i:1\leq i\leq 2^b-1)\cdot y_0,
\end{equation*}
while the subfield $K(U_i,V_i:1\leq i\leq 2^b-1)$ is invariant by
the action of $\widetilde G$. Thus $K(x_i,y_i:0\leq i\leq
2^b-1)^{\widetilde G}=K(U_i,V_i:1\leq i\leq 2^b-1)^{\widetilde
G}(u,v)$ for some $u,v$ such that $\sigma(v)=\tau(v)=\rho(v)=v$ and
$\sigma(u)=\tau(u)=\rho(u)=u$. We have now {\allowdisplaybreaks
\begin{align*}
\sigma\ :\ &U_i\mapsto \xi^{k^{i-1}} U_i,~ V_i\mapsto\xi^{k^i-k^{i-1}}V_i,~ v\mapsto v\\
\tau\ :\ &U_1\mapsto U_2\mapsto\cdots\mapsto U_{2^b-1}\mapsto (U_1U_2\cdots U_{2^b-1})^{-1},\\
&V_1\mapsto V_2\mapsto\cdots\mapsto V_{2^b-1}\mapsto (V_1V_2\cdots V_{2^b-1})^{-1},\\
&~ v\mapsto v,\\
\rho\ :\ &U_i\mapsto U_i,~ V_i\mapsto V_i,~ v\mapsto v
\end{align*}}
for $1\leq i\leq 2^b-1$. From Theorem \ref{t2.2} it follows that if
$K(U_i,V_i:1\leq i\leq 2^b-1)^{\widetilde G}(v)$ is rational over
$K$, so is $K(x_i,y_i:0\leq i\leq 2^b-1)^{\widetilde G}$ over $K$.

Since $\rho$ acts trivially on $K(U_i,V_i:1\leq i\leq 2^b-1)$, we
find that $K(U_i,V_i:1\leq i\leq 2^b-1)^{\widetilde
G}=K(U_i,V_i:1\leq i\leq 2^b-1)^{\langle\sigma,\tau\rangle}$.

For $1\leq i\leq 2^b-1$, define $v_i=V_i/U_i^{k-1}$. Let $\eta\in K$
be a primitive $2^{t+1}$-th root of unity such that
$\xi=\eta^{-2+2^r}$. Then $\xi^{k^{i-1}}=\eta^{k^i-k^{i-1}}$ for any
$i$. We have now the actions {\allowdisplaybreaks
\begin{align*}
\sigma\ :\ &U_i\mapsto \eta^{k^i-k^{i-1}} U_i,~ v_i\mapsto v_i,~ v\mapsto v\\
\tag{5.4} \tau\ :\ &U_1\mapsto U_2\mapsto\cdots\mapsto U_{2^b-1}\mapsto (U_1U_2\cdots U_{2^b-1})^{-1},\\
&v_1\mapsto v_2\mapsto\cdots\mapsto v_{2^b-1}\mapsto (v_1v_2\cdots v_{2^b-1})^{-1},\\
&v\mapsto v,
\end{align*}}
for $1\leq i\leq 2^b-1$. Compare Formula (5.4) with Formula (5.3).
They look almost the same. Apply the proof of Case 3.

{\bf Case 6.} $\widetilde G=G_6$. Similarly to Case 5, we may assume
that $a=t$. The proof is almost the same as in Case 5. Here only the
action of $\tau$ is changed. As we already saw in Case 2, this issue
is solved with a proper adjustment involving certain root of unity
which is in $K$.

{\bf Case 7.} $\widetilde G=G_7$. From the relation
$\tau^{-1}\sigma^{2^a}\tau=\sigma^{-2^a+2^{a+r}}=\sigma^{2^a}$ it
follows that the order of $\sigma$ is $2^{a+1}$. Therefore,
$\alpha=t-1$. Note that the subgroup $H=\langle\sigma,\rho\rangle$
is abelian and has an order $p^{a+t}$. Put
$\rho_1=\sigma^{2^{a+1-t}}\rho^{-1}$. Then
$H\cong\langle\sigma\rangle\times\langle\rho_1\rangle$, where
$\sigma^{2^{a+1}}=\rho_1^{2^{t-1}}=1$. We have also
$\rho=\sigma^{2^{a+1-t}}\rho_1^{-1}$.

Define $X_1,X_2\in V^*$ by
\begin{equation*}
X_1=\sum_{i=0}^{2^{a+1}-1}x(\sigma^i),~
X_2=\sum_{i=0}^{2^{t-1}-1}x(\rho_1^i).
\end{equation*}
Note that $\sigma\cdot X_1=X_1$ and $\rho_1\cdot X_2=X_2$.

Let $\zeta_1$ be a primitive $2^{a+1}$-th root of unity;
$\xi=\zeta^{2^{a+1-t}}$, a primitive $2^t$-th root of unity; and
$\mu=\zeta^{2^{a+1-t+1}}$, a primitive $2^{t-1}$-th root of unity.
Define $Y_1,Y_2\in V^*$ by
\begin{equation*}
Y_1=\sum_{i=0}^{2^{t-1}-1}\mu^{-i}\rho_1^i\cdot X_1,~
Y_2=\sum_{i=0}^{2^{a+1}-1}\zeta_1^{-i}\sigma^i\cdot X_2.
\end{equation*}

It follows that {\allowdisplaybreaks
\begin{align*}
\sigma\ :\ &Y_1\mapsto Y_1,Y_2\mapsto\zeta_1 Y_2,\\
\rho\ :\ &Y_1\mapsto\mu^{-1} Y_1,Y_2\mapsto\xi Y_2.
\end{align*}}
Thus $K\cdot Y_1+K\cdot Y_2$ is a representation space of the
subgroup $H$.

Define $x_i=\tau^i\cdot Y_1,y_i=\tau^i\cdot Y_2$ for $0\leq i\leq
2^b-1$. We have now {\allowdisplaybreaks \begin{align*}
\sigma\ :\ &x_i\mapsto \mu^{-w_i} x_i,~ y_i\mapsto\xi^{w_i}\zeta_1^{k^i}y_i\\
\tau\ :\ &x_0\mapsto x_1\mapsto\cdots\mapsto x_{2^b-1}\mapsto x_0,\\
&y_0\mapsto y_1\mapsto\cdots\mapsto y_{2^b-1}\mapsto y_0,\\
\rho\ :\ &x_i\mapsto\mu^{-1} x_i,~ y_i\mapsto\xi y_i,
\end{align*}}
for $0\leq i\leq 2^b-1$. Computations show that for $a=t$ we have
$\sigma^{2^{t-1}}(y_i)=\rho^{2^{t-2}}(y_i)$ for any $i$, which means
that $Y=\bigoplus_{0\leq i\leq 2^b-1}K\cdot y_i$ is not a faithful
$\widetilde G$-subspace of $V^*$.

Define $z_i=x_iy_i$ for $1\leq i\leq 2^b-1$. It follows that
{\allowdisplaybreaks
\begin{align*}
\sigma\ :\ &x_i\mapsto \mu^{-w_i} x_i,~ z_i\mapsto(\mu^{-1}\xi)^{w_i}\zeta_1^{k^i}z_i\\
\tau\ :\ &x_0\mapsto x_1\mapsto\cdots\mapsto x_{2^b-1}\mapsto x_0,\\
&z_0\mapsto z_1\mapsto\cdots\mapsto z_{2^b-1}\mapsto z_0,\\
\rho\ :\ &x_i\mapsto\mu^{-1} x_i,~ z_i\mapsto\mu^{-1}\xi z_i,
\end{align*}}
for $0\leq i\leq 2^b-1$. Now for $a=t$ we have
$\sigma^{2^{t-1}}(z_i)=-\rho^{2^{t-2}}(z_i)$, and for any other
value of $a>t$ we have
$\sigma^{2^{a-1}}(z_i)\ne\rho^{2^{t-2}}(z_i)$. It follows that
$Y=\bigoplus_{0\leq i\leq 2^b-1}K\cdot z_i$ is a faithful
$\widetilde G$-subspace of $V^*$. Thus, by Theorem \ref{t2.1}, it
suffices to show that $K(z_i:0\leq i\leq 2^b-1)^{\widetilde G}$ is
rational over $K$.

For $1\leq i\leq 2^b-1$ define $U_i=z_i/z_{i-1}$. We have now
{\allowdisplaybreaks
\begin{align*}
\sigma\ :\ &U_i\mapsto (\mu^{-1}\xi)^{k^{i-1}}\zeta_1^{k^i-k^{i-1}} U_i,\\
\tag{5.5} \tau\ :\ &U_1\mapsto U_2\mapsto\cdots\mapsto
U_{2^b-1}\mapsto
(U_1U_2\cdots U_{2^b-1})^{-1},\\
\rho\ :\ &U_i\mapsto U_i,
\end{align*}}
for $1\leq i\leq 2^b-1$.

Let $\eta$ be a primitive $2^{t+1}$-th root of unity such that
$\mu^{-1}\xi=\eta^{k-1}$. Whence
$(\mu^{-1}\xi)^{k^{i-1}}=\eta^{k^i-k^{i-1}}$.

Compare Formula (5.5) with Formula (5.3). They look almost the same.
Apply the proof of Case 3.

{\bf Cases 8-16.} $\widetilde G=G_i$  for $8\leq i\leq 16$. It is
easily seen that we need to make only minor changes, as we did in
Cases 2 and 8, in order to ensure the proper action of $\tau$.

{\bf Step II.} Consider the general presentation of $\widetilde G$.
According to Proposition \ref{p3.2}, we may assume that $\widetilde
G$ has the following presentation:
\begin{equation*}
\widetilde
G=\langle\sigma,\tau,\rho:\sigma^{p^a}=\rho^{sp^\alpha},\tau^{p^b}=\sigma^{mp^c}\rho^{p^\beta},\rho^{p^t}=1,\tau^{-1}\sigma\tau=\sigma^k\rho,
\rho - \text{central}\rangle,
\end{equation*}
where $s,m$ are positive integers, $1\leq s<p^t,0 \leq m<p^a$,
$\gcd(sm,p)=1,0\leq\alpha,\beta\leq t$ and $k=\varepsilon+p^r$.

Here we have again $16$ cases, which correspond to those in Step I.
Since they all can be treated in an unified way, we will consider
only Case 3.

Let $c=\beta=t$, i.e., $\tau^{p^b}=1$. The subgroup
$H=\langle\sigma,\rho\rangle$ is abelian and has an order $p^{a+t}$.
Put $\rho_1=\sigma^{p^{a-\alpha}}\rho^{-s}$. Then
$H\cong\langle\sigma\rangle\times\langle\rho_1\rangle$, where
$\sigma^{p^{a+t-\alpha}}=\rho_1^{p^\alpha}=1$. Let $n$ be an integer
such that $ns\equiv 1\pmod{p^t}$. We have now
$\rho=(\rho^s)^n=\sigma^{np^{a-\alpha}}\rho_1^{-n}$.

Define $X_1,X_2\in V^*$ by
\begin{equation*}
X_1=\sum_ix(\rho_1^i),~ X_2=\sum_ix(\sigma^i).
\end{equation*}
Note that $\sigma\cdot X_2=X_2,\rho_1\cdot X_1=X_1$.

Let $\zeta_1\in K$ be a primitive $p^{a+t-\alpha}$-th root of unity.
Put $\xi=\zeta_1^{p^{a-\alpha}}$, a primitive $p^t$-th root of
unity. Let $\zeta_2$ be any primitive $p^\alpha$-th root of unity.
(We will specify $\zeta_2$ a bit later.)

Define $Y_1,Y_2,Y_3\in V^*$ by
\begin{equation*}
Y_1=\sum_{i=0}^{p^{a+t-\alpha}-1}\zeta_1^{-i}\sigma^i\cdot X_1,~
Y_2=\sum_{i=0}^{p^\alpha-1}\zeta_2^{-i}\rho_1^i\cdot X_3.
\end{equation*}

It follows that {\allowdisplaybreaks \begin{align*}
\sigma\ :\ &Y_1\mapsto\zeta_1 Y_1,~ Y_2\mapsto Y_2,\\
\rho_1\ :\ &Y_1\mapsto Y_1,~ Y_2\mapsto\zeta_2 Y_2,\\
\rho\ :\ &Y_1\mapsto\xi^n Y_1,~ Y_2\mapsto\zeta_2^{-n} Y_2.
\end{align*}}
Thus $K\cdot Y_1+K\cdot Y_2$ is a representation space of the
subgroup $H$.

Define $x_i=\tau^i\cdot Y_1,y_i=\tau^i\cdot Y_2$ for $0\leq i\leq
p^b-1$. We have now {\allowdisplaybreaks
\begin{align*}
\sigma\ :\ &x_i\mapsto \zeta_1^{k^i}\xi^{nw_i} x_i,~ y_i\mapsto\zeta_2^{-nw_i}y_i,\\
\tau\ :\ &x_0\mapsto x_1\mapsto\cdots\mapsto x_{p^b-1}\mapsto x_0,\\
&y_0\mapsto y_1\mapsto\cdots\mapsto y_{p^b-1}\mapsto y_0,\\
\rho\ :\ &x_i\mapsto\xi^n x_i,~ y_i\mapsto\zeta_2^{-n} y_i.
\end{align*}}
for $0\leq i\leq p^b-1$. We find that $Y=(\bigoplus_{0\leq i\leq
p^b-1}K\cdot x_i)\oplus(\bigoplus_{0\leq i\leq p^b-1}K\cdot y_i)$ is
a faithful $\widetilde G$-subspace of $V^*$. Thus, by Theorem
\ref{t2.1}, it suffices to show that $K(x_i,y_i:0\leq i\leq
p^b-1)^{\widetilde G}$ is rational over $K$.

For $1\leq i\leq p^b-1$, define $U_i=x_i/x_{i-1}$ and
$V_i=y_i/y_{i-1}$. Thus $K(x_i,y_i:0\leq i\leq
p^b-1)=K(x_0,y_0,U_i,V_i:1\leq i\leq p^b-1)$ and for every
$g\in\widetilde G$
\begin{equation*}
g\cdot x_0\in K(U_i,V_i)\cdot x_0,~ g\cdot y_0\in K(U_i,V_i)\cdot
y_0,
\end{equation*}
while the subfield $K(U_i,V_i:1\leq i\leq p^b-1)$ is invariant by
the action of $\widetilde G$, i.e., {\allowdisplaybreaks
\begin{align*}
\sigma\ :\ &U_i\mapsto \zeta_1^{k^i-k^{i-1}}\xi^{nk^{i-1}} U_i,~ V_i\mapsto\zeta_2^{-nk^{i-1}} V_i,\\
\tau\ :\ &U_1\mapsto U_2\mapsto\cdots\mapsto U_{p^b-1}\mapsto (U_1\cdots U_{p^b-1})^{-1},\\
&V_1\mapsto V_2\mapsto\cdots\mapsto V_{p^b-1}\mapsto (V_1\cdots V_{p^b-1})^{-1},\\
\rho\ :\ &U_i\mapsto U_i,~ V_i\mapsto V_i.
\end{align*}}
for $1\leq i\leq p^b-1$. From Theorem \ref{t2.2} it follows that if
$K(U_i,V_i:1\leq i\leq p^b-1)^{\widetilde G}$ is rational over $K$,
so is $K(x_i,y_i:0\leq i\leq p^b-1)^{\widetilde G}$ over $K$.

Since $\rho$ acts trivially on $K(U_i,V_i)$, we find that
$K(U_i,V_i)^{\widetilde G}=K(U_i,V_i)^{\langle\sigma,\tau\rangle}$.

{\bf Subcase 3.a.} Let $a-\alpha\leq r$. Thus we can write
$r=a-\alpha+r_1$ for some $r_1\geq 0$. Therefore,
$\zeta_1^{k^{i-1}(k-1)}=\xi^{k^{i-1}p^{r_1}}$ for all $i$. Define
$\zeta_2=\xi^{(1+sp^{r_1})p^{t-\alpha}}$, a primitive $p^\alpha$-th
root of unity.

Define $v_i=U_i^{p^{t-\alpha}}V_i$. Since
$\xi^{(n+p^{r_1})p^{t-\alpha}}=\xi^{n(1+sp^{r_1})p^{t-\alpha}}=\zeta_2^n$,
we have {\allowdisplaybreaks
\begin{align*}
\sigma\ :\ &U_i\mapsto \xi^{k^{i-1}n(1+sp^{r_1})} U_i,~ v_i\mapsto v_i,\\
\tag{5.6} \tau\ :\ &U_1\mapsto U_2\mapsto\cdots\mapsto U_{p^b-1}\mapsto (U_1\cdots U_{p^b-1})^{-1},\\
&v_1\mapsto v_2\mapsto\cdots\mapsto v_{p^b-1}\mapsto (v_1\cdots
v_{p^b-1})^{-1}
\end{align*}}
for $1\leq i\leq p^b-1$. Compare Formula (5.6) with Formula (5.1).
They look almost the same. Apply the proof of Case 1.

{\bf Subcase 3.b.} Let $a-\alpha>r$. Let $\eta$ be a primitive
$p^{r+t}$-th root of unity such that $\xi^n=\eta^{p^r}$. Whence
$\xi^{nk^{i-1}}=\eta^{k^i-k^{i-1}}$. Since $a+t-\alpha>r+t$, we get
that $\zeta_1\eta$ is a primitive $p^{a+t-\alpha}$-th root of unity.
Put $\mu=(\zeta_1\eta)^{p^{a+t-\alpha-(\alpha+r)}}$, a primitive
$p^{\alpha+r}$-th root of unity, where $\alpha+r<a\leq a+t-\alpha$.
Now, put $\zeta_2=\mu^{sp^r}$, a primitive $p^\alpha$-th root of
unity (recall that $ns\equiv 1\pmod{p^t}$). Thus
$\zeta_2^{nk^{i-1}}=\mu^{k^i-k^{i-1}}$ for any $i$.

Define $v_i=U_i^{p^{a+t-\alpha-(\alpha+r)}}V_i$ for $1\leq i\leq
p^b-1$. It follows that {\allowdisplaybreaks
\begin{align*}
\sigma\ :\ &U_i\mapsto (\zeta_1\eta)^{k^i-k^{i-1}} U_i,~ v_i\mapsto v_i,\\
\tag{5.7} \tau\ :\ &U_1\mapsto U_2\mapsto\cdots\mapsto U_{p^b-1}\mapsto (U_1\cdots U_{p^b-1})^{-1},\\
&v_1\mapsto v_2\mapsto\cdots\mapsto v_{p^b-1}\mapsto (v_1\cdots
v_{p^b-1})^{-1}
\end{align*}}
for $1\leq i\leq p^b-1$. Compare Formula (5.7) with Formula (5.3).
They are the same. Apply the proof of Case 3, Step I.

\section{Proof of Corollary \ref{c1.6}}
\label{6}

Let $C=C_{p^{\alpha_1}}\times\cdots\times C_{p^{\alpha_s}}\leq
Z(\widetilde G)$. Denote again by $\sigma$ and $\tau$ the preimages
of the generators of $G$ and by $\rho_1,\dots,\rho_s$ the generators
of $C$, i.e., $\rho_i^{p^{\alpha_i}}=1$. Then
$[\sigma,\tau]=\rho_1^{\beta_1}\cdots\rho_s^{\beta_s}$ for
$\beta_i\geq 0$, and we can assume for abuse of notation that
$p^t=\ord(\rho_1^{\beta_1})=\max\{\ord(\rho_1^{\beta_1}),\dots,\ord(\rho_s^{\beta_s})\}$.
It follows that $\widetilde
G'\cap\langle\rho_2,\dots,\rho_s\rangle=\{1\}$. Now we can apply
Theorem \ref{t2.8} reducing the rationality problem of $K(\widetilde
G)$ over $K$ to the rationality problem of $K(\widetilde
G/\langle\rho_2,\dots,\rho_k\rangle)$ over $K$. Put $\widetilde
G_1=\widetilde G/\langle\rho_2,\dots,\rho_k\rangle$ and let us find
an upper bound for the exponent of $\widetilde G_1/\widetilde G_1'$.
We may assume that $\widetilde G_1$ has the presentation
\begin{equation*}
\widetilde
G_1=\langle\sigma,\tau,\rho:\sigma^{p^a}=\rho^{sp^\alpha},\tau^{p^b}=\sigma^{mp^c}\rho^{p^\beta},\rho^{p^t}=1,\tau^{-1}\sigma\tau=\sigma^k\rho,
\rho - \text{central}\rangle,
\end{equation*}
where $s,m$ are positive integers, $1\leq s<p^t,0 \leq m<p^a$,
$\gcd(sm,p)=1,0\leq\alpha,\beta\leq t$ and $k=\varepsilon+p^r$.

Computations show that we will obtain an upper bound for the
exponent of $\widetilde G_1$ in the following situation:
$\sigma^{p^a}=\rho^{p^\alpha},\tau^{p^b}=\sigma^{p^c}\rho^{p^\beta}$
for $\alpha<t,c<a,\beta<t$. We have
$\tau^{p^{b+a-c}}=\rho^{p^\alpha+p^{\beta+a-c}}$. There are several
possibilities for the order of $\tau$.

{\bf Case 1.} Let $\beta+a-c\geq t$. Then
$\ord(\tau)=p^{b+a+t-c-\alpha}$.

{\bf Case 2.} Let $\beta+a-c<t$ and $\alpha<\beta+a-c$. Here again
$\ord(\tau)=p^{b+a+t-c-\alpha}$.

{\bf Case 3.} Let $\beta+a-c<t$ and $\alpha\geq\beta+a-c$. We have
$\ord(\tau)=p^{b+t-\beta}$.

Note also that $\ord(\sigma)=p^{a+t-\alpha}\leq p^{a+t}\leq
p^{a+b+t-c}$.

Therefore, $\exp(\widetilde G_1/\widetilde G_1')\leq\exp(\widetilde
G_1)\leq p^{a+b+t-c}$.

\section{Proof of Theorem \ref{t1.7}}
\label{7}

The case (i) follows from Theorem \ref{t2.3}.

(ii) We will divide the proof into two steps.

{\bf Step I.} Assume that $\widetilde G$ has the following
presentation:
\begin{equation*}
\widetilde
G=\langle\sigma,\tau,\rho:\sigma^{p^a}=\rho^{p^\alpha},\tau^{p^b}=\sigma^{p^c}\rho^{p^\beta},\rho^{p^t}=1,\tau^{-1}\sigma\tau=\sigma^k\rho,
\rho - \text{central}\rangle,
\end{equation*}
where $a\geq t, b\geq t, r\geq t, 0\leq\alpha,\beta\leq t$ and
$k=\varepsilon+p^r$. Note that for any $i$ we have the relation
$\tau^{-i}\sigma\tau^{i}=\sigma^{k^i}\rho^i$, since $\rho^k=\rho$.

{\bf Case 1.} $\widetilde G=G_1$, where $G_1$ is the group in
Section \ref{3}. From Proposition \ref{p3.3} it follows that
$\widetilde G'$ is cyclic and is generated by $\sigma^{p^r}\rho$.

{\bf Subcase 1.a.} Let $a\geq r+t$. Then
$(\sigma^{p^r}\rho)^{p^{a-r}}=\rho^{p^{a-r}}=1$, so $\widetilde
G'\cap\langle\rho\rangle=1$. Theorem \ref{t2.8} then implies that we
can reduce the rationality problem of $K(\widetilde G)$ to
$K(\widetilde G/\langle\rho\rangle)$ over $K$, where $\widetilde
G/\langle\rho\rangle$ clearly is a metacyclic $p$-group.

{\bf Subcase 1.b.} Let $a<r+t$.  Then
$\tau^{-p^t}\sigma\tau^{p^t}=\sigma^{k^{p^t}}=\sigma$, since
$(1+p^r)^{p^t}=1+A\cdot p^{r+t}$ for some integer $A$. This means
that $\tau^{p^t}$ is in the center of $\widetilde G$. Clearly,
$\widetilde G'\cap\langle\tau^{p^t}\rangle=1$, so we can apply again
theorem \ref{t2.8} reducing the rationality problem of $K(\widetilde
G)$ to $K(\widetilde G/\langle\tau^{p^t}\rangle)$ over $K$. In other
words, we may assume that $\tau^{p^t}=1$.

Let $V$ be a $K$-vector space whose dual space $V^*$ is defined as
$V^*=\bigoplus_{g\in\widetilde G}K\cdot x(g)$ where $\widetilde G$
acts on $V^*$ by $h\cdot x(g)=x(hg)$ for any $h,g\in\widetilde G$.
Thus $K(V)^{\widetilde G}=K(x(g):g\in\widetilde G)^{\widetilde
G}=K(\widetilde G)$.

Define $X_1,X_2\in V^*$ by
\begin{equation*}
X_1=\sum_{i=0}^{p^a-1}x(\sigma^i),~ X_2=\sum_{i=0}^{p^t-1}x(\rho^i).
\end{equation*}
Note that $\sigma\cdot X_1=X_1$ and $\rho\cdot X_2=X_2$.

Let $\zeta=\zeta_{p^a}\in K$ be a primitive $p^a$-th root of unity
and define $\xi=\zeta^{p^{a-t}}$. Thus $\xi$ is a primitive $p^t$-th
root of unity. Define $Y_1,Y_2\in V^*$ by
\begin{equation*}
Y_1=\sum_{i=0}^{p^t-1}\xi^{-i}\rho^i\cdot X_1,~
Y_2=\sum_{i=0}^{p^a-1}\zeta^{-i}\sigma^i\cdot X_2.
\end{equation*}

It follows that {\allowdisplaybreaks \begin{align*}
\sigma\ :\ &Y_1\mapsto Y_1,Y_2\mapsto\zeta Y_2,\\
\rho\ :\ &Y_1\mapsto\xi Y_1,Y_2\mapsto Y_2.
\end{align*}}
Thus $K\cdot Y_1+K\cdot Y_2$ is a representation space of the
subgroup $\langle\sigma,\rho\rangle$.

Define $x_i=\tau^i\cdot Y_1,y_i=\tau^i\cdot Y_2$ for $0\leq i\leq
p^t-1$. We have now {\allowdisplaybreaks \begin{align*}
\sigma\ :\ &x_i\mapsto \xi^i x_i,~ y_i\mapsto\zeta^{k^i}y_i\\
\tau\ :\ &x_0\mapsto x_1\mapsto\cdots\mapsto x_{p^t-1}\mapsto x_0,\\
&y_0\mapsto y_1\mapsto\cdots\mapsto y_{p^t-1}\mapsto y_0,\\
\rho\ :\ &x_i\mapsto\xi x_i,~ y_i\mapsto y_i,
\end{align*}}
for $0\leq i\leq p^t-1$.

We find that $Y=(\bigoplus_{0\leq i\leq p^t-1}K\cdot
x_i)\oplus(\bigoplus_{0\leq i\leq p^t-1}K\cdot y_i)$ is a faithful
$\widetilde G$-subspace of $V^*$. Thus, by Theorem \ref{t2.1}, it
suffices to show that $K(x_i,y_i:0\leq i\leq p^t-1)^{\widetilde G}$
is rational over $K$.

For $1\leq i\leq p^t-1$, define $u_i=x_i/x_{i-1}$ and
$v_i=y_i/y_{i-1}$. Thus $K(x_i,y_i:0\leq i\leq
p^t-1)=K(x_0,y_0,u_i,v_i:1\leq i\leq p^t-1)$ and for every
$g\in\widetilde G$
\begin{equation*}
g\cdot x_0\in K(u_i,v_i:1\leq i\leq p^t-1)\cdot x_0,~ g\cdot y_0\in
K(u_i,v_i:1\leq i\leq p^t-1)\cdot y_0,
\end{equation*}
while the subfield $K(u_i,v_i:1\leq i\leq p^t-1)$ is invariant by
the action of $\widetilde G$, i.e., {\allowdisplaybreaks
\begin{align*}
\sigma\ :\ &u_i\mapsto \xi u_i,~ v_i\mapsto\zeta^{k^i-k^{i-1}}v_i\\
\tau\ :\ &u_1\mapsto u_2\mapsto\cdots\mapsto u_{p^t-1}\mapsto (u_1u_2\cdots u_{p^t-1})^{-1},\\
&v_1\mapsto v_2\mapsto\cdots\mapsto v_{p^t-1}\mapsto (v_1v_2\cdots v_{p^t-1})^{-1},\\
\rho\ :\ &u_i\mapsto u_i,~ v_i\mapsto v_i,
\end{align*}}
for $1\leq i\leq p^t-1$. From Theorem \ref{t2.2} it follows that if
$K(u_i,v_i:1\leq i\leq p^t-1)^{\widetilde G}$ is rational over $K$,
so is $K(x_i,y_i:0\leq i\leq p^t-1)^{\widetilde G}$ over $K$.

Since $\rho$ acts trivially on $K(u_i,v_i:1\leq i\leq p^t-1)$, we
find that $K(u_i,v_i:1\leq i\leq p^t-1)^{\widetilde
G}=K(u_i,v_i:1\leq i\leq p^t-1)^{\langle\sigma,\tau\rangle}$.

Recall that $a<r+t$. So we can write $r=a-t+r_1$ for some $r_1\geq
1$. Since $\xi^k=\xi$, we have
$\zeta^{k^{i-1}(k-1)}=\zeta^{k^{i-1}p^{a-t+r_1}}=\xi^{p^{r_1}}$ for
$1\leq i\leq p^t-1$.

Define $w_i=v_i/(u_i)^{p^{r_1}}$. We have now that $K(u_i,v_i:1\leq
i\leq p^t-1)=K(u_i,w_i:1\leq i\leq p^t-1)$ and {\allowdisplaybreaks
\begin{align*}
\sigma\ :\ &u_i\mapsto \xi u_i,~ w_i\mapsto w_i\\
\tau\ :\ &u_1\mapsto u_2\mapsto\cdots\mapsto u_{p^t-1}\mapsto (u_1u_2\cdots u_{p^t-1})^{-1},\\
&w_1\mapsto w_2\mapsto\cdots\mapsto w_{p^t-1}\mapsto (w_1w_2\cdots
w_{p^t-1})^{-1}.
\end{align*}}
According to Lemma \ref{l2.7}, we can linearize the action of $\tau$
on $w_1,\dots,w_{p^t-1}$.

Write $L=K(u_i:1\leq i\leq p^t-1)$ and consider $L(w_i:0\leq i\leq
p^t-1)^{\langle\sigma,\tau\rangle}$. Note that the group
$\langle\sigma,\tau\rangle$ acts on the field $L(w_i)$ as
$\langle\sigma,\tau\rangle/\langle\sigma^{p^t}\rangle$ and is
faithful on $L$. Thus we may apply Theorem \ref{t2.1} to
$L(w_i:1\leq i\leq p^t-1)^{\langle\sigma,\tau\rangle}$. It remains
to show that $L^{\langle\sigma,\tau\rangle}$ is rational over $K$.

Define $z_1=u_1^{p^t},z_i=u_i/u_{i-1}$ for $2\leq i\leq p^t-1$. Then
$L^{\langle\sigma\rangle}=K(z_i:1\leq i\leq p^t-1)$ and the action
of $\tau$ is given by {\allowdisplaybreaks
\begin{align*}
\tau\ :\ &z_1\mapsto z_1z_2^{p^t},\\
&z_2\mapsto z_3\mapsto\cdots\mapsto z_{p^t-1}\mapsto
(z_1z_2^{p^t-1}z_3^{p^t-2}\cdots z_{p^t-1}^2)^{-1}\mapsto
z_1z_2^{p^t-2}z_3^{p^t-3}\cdots z_{p-2}^2z_{p^t-1}\mapsto z_2.
\end{align*}}
Define $s_1=z_2,s_i=\tau^{i-1}\cdot z_2$ for $2\leq i\leq p^t-1$.
Then $K(z_i:1\leq i\leq p^t-1)=K(s_i:1\leq i\leq p^t-1)$ and
\begin{align*}
\tau\ :\ &s_1\mapsto s_2\mapsto\cdots\mapsto s_{p^t-1}\mapsto
(s_1s_2\cdots s_{p^t-1})^{-1}.
\end{align*}
The action of $\tau$ can be linearized according to Lemma
\ref{l2.7}. Thus $K(s_i:1\leq i\leq p^t-1)^{\langle\tau\rangle}$ is
rational over $K$ by Theorem \ref{t1.1}.

{\bf Case 2.} $\widetilde G=G_2$. We can apply the same argument of
Subcase 1.a, so we will assume that $t+r>a$. Note that
$\tau^{-p^t}\sigma\tau^{p^t}=\sigma^{(1+p^r)^{p^t}}\rho^{p^t}=\sigma$,
since $t+r>a$. Therefore, $\tau^{p^t}$ is in the center of
$\widetilde G$ and the group
$H=\langle\sigma,\rho,\tau^{p^t}\rangle$ is abelian of order
$p^{a+b}$. There are two possibilities for the decomposition of $H$
as a direct product of cyclic groups.

{\bf Subcase 2.a.} $b-\beta\geq t$. Define
$\rho_1=\tau^{p^{b-\beta}}\rho^{-1},\rho_2=\tau^{p^t}$. Then $H$ is
isomorphic to the direct product
$\langle\sigma\rangle\times\langle\rho_1\rangle\times\langle\rho_2\rangle$,
where $\sigma^{p^a}=\rho_1^{p^\beta}=\rho_2^{p^{b-\beta}}=1$. We
have $\rho=\rho_1^{-1}\rho_2^{p^{b-\beta-t}}$.

Define $X_1,X_2,X_3\in V^*$ by
\begin{equation*}
X_1=\sum_{i,j}x(\rho_1^i\rho_2^j),~
X_2=\sum_{i,j}x(\sigma^i\rho_1^j),~
X_3=\sum_{i,j}x(\sigma^i\rho_2^j).
\end{equation*}
Note that $\sigma\cdot X_2=X_2,\sigma\cdot X_3=X_3,\rho_1\cdot
X_1=X_1,\rho_1\cdot X_2=X_2,\rho_2\cdot X_1=X_1$ and $\rho_2\cdot
X_3=X_3$.

Let $\zeta=\zeta_{p^a}\in K$ be a primitive $p^a$-th root of unity.
Define $\zeta_1=\zeta^{p^{a-\beta}}$, a primitive $p^\beta$-th root
of unity; $\zeta_2=\zeta^{p^{a-b+\beta}}$, a primitive
$p^{b-\beta}$-th root of unity; and $\xi=\zeta^{p^{a-t}}$, a
primitive $p^t$-th root of unity.

Define $Y_1,Y_2,Y_3\in V^*$ by
\begin{equation*}
Y_1=\sum_{i=0}^{p^a-1}\zeta^{-i}\sigma^i\cdot X_1,~
Y_2=\sum_{i=0}^{p^\beta-1}\zeta_1^{-i}\rho_1^i\cdot X_3,~
Y_3=\sum_{i=0}^{p^{b-\beta}-1}\zeta_2^{-i}\rho_2^i\cdot X_2.
\end{equation*}

It follows that {\allowdisplaybreaks \begin{align*}
\sigma\ :\ &Y_1\mapsto\zeta Y_1,~ Y_2\mapsto Y_2,~ Y_3\mapsto Y_3,\\
\rho_1\ :\ &Y_1\mapsto Y_1,~ Y_2\mapsto\zeta_1 Y_2,~ Y_3\mapsto Y_3,\\
\rho_2\ :\ &Y_1\mapsto Y_1,~ Y_2\mapsto Y_2,~ Y_3\mapsto\zeta_2 Y_3,\\
\rho\ :\ &Y_1\mapsto Y_1,~ Y_2\mapsto\zeta_1^{-1} Y_2,~
Y_3\mapsto\xi Y_3.
\end{align*}}
Thus $K\cdot Y_1+K\cdot Y_2+K\cdot Y_3$ is a representation space of
the subgroup $H$.

Define $x_i=\tau^i\cdot Y_1,y_i=\tau^i\cdot Y_2,z_i=\tau^i\cdot Y_3$
for $0\leq i\leq p^t-1$. We have now {\allowdisplaybreaks
\begin{align*}
\sigma\ :\ &x_i\mapsto \zeta^{k^i} x_i,~ y_i\mapsto\zeta_1^{-i}y_i,~ z_i\mapsto\xi^iz_i,\\
\tau\ :\ &x_0\mapsto x_1\mapsto\cdots\mapsto x_{p^t-1}\mapsto x_0,\\
&y_0\mapsto y_1\mapsto\cdots\mapsto y_{p^t-1}\mapsto y_0,\\
&z_0\mapsto z_1\mapsto\cdots\mapsto z_{p^t-1}\mapsto\zeta_2 z_0,\\
\rho\ :\ &x_i\mapsto x_i,~ y_i\mapsto\zeta_1^{-1} y_i,~
z_i\mapsto\xi z_i,
\end{align*}}
for $0\leq i\leq p^t-1$. We find that $Y=(\bigoplus_{0\leq i\leq
p^t-1}K\cdot x_i)\oplus(\bigoplus_{0\leq i\leq p^t-1}K\cdot
y_i)\oplus(\bigoplus_{0\leq i\leq p^t-1}K\cdot z_i)$ is a faithful
$\widetilde G$-subspace of $V^*$. Thus, by Theorem \ref{t2.1}, it
suffices to show that $K(x_i,y_i,z_i:0\leq i\leq p^t-1)^{\widetilde
G}$ is rational over $K$.

Note that $\zeta_1=\xi^{p^{t-\beta}}$. For $1\leq i\leq p^t-1$,
define $U_i=x_i/x_{i-1}$ and $W_i=z_i/z_{i-1}$. For $0\leq i\leq
p^t-1$, define $V_i=y_iz_i^{p^{t-\beta}}$. Thus $K(x_i,y_i,z_i:0\leq
i\leq p^t-1)=K(x_0,z_0,U_i,V_j,W_i:1\leq i\leq p^t-1,0\leq j\leq
p^t-1)$ and for every $g\in\widetilde G$
\begin{equation*}
g\cdot x_0\in K(U_i,V_j,W_i)\cdot x_0,~ g\cdot z_0\in
K(U_i,V_j,W_i)\cdot z_0,
\end{equation*}
while the subfield $K(U_i,V_j,W_i:1\leq i\leq p^t-1,0\leq j\leq
p^t-1)$ is invariant by the action of $\widetilde G$, i.e.,
{\allowdisplaybreaks
\begin{align*}
\sigma\ :\ &U_i\mapsto \zeta^{k^i-k^{i-1}} U_i,~ V_j\mapsto V_j,~ W_i\mapsto\xi W_i,\\
\tau\ :\ &U_1\mapsto U_2\mapsto\cdots\mapsto U_{p^t-1}\mapsto (U_1\cdots U_{p^t-1})^{-1},\\
\tag{7.1} &V_0\mapsto V_1\mapsto\cdots\mapsto V_{p^t-1}\mapsto\zeta_2^{p^{t-\beta}} V_0,\\
&W_1\mapsto W_2\mapsto\cdots\mapsto W_{p^t-1}\mapsto\zeta_2 (W_1\cdots W_{p^t-1})^{-1},\\
\rho\ :\ &U_i\mapsto U_i,~ V_j\mapsto V_j,~ W_i\mapsto W_i,
\end{align*}}
for $1\leq i\leq p^t-1,0\leq j\leq p^t-1$. From Theorem \ref{t2.2}
it follows that if $K(U_i,V_j,W_i:1\leq i\leq p^t-1,0\leq j\leq
p^t-1)^{\widetilde G}$ is rational over $K$, so is
$K(x_i,y_i,z_i:0\leq i\leq p^t-1)^{\widetilde G}$ over $K$.

Since $\rho$ acts trivially on $K(U_i,V_j,W_i)$, we find that
$K(U_i,V_j,W_i)^{\widetilde
G}=K(U_i,V_j,W_i)^{\langle\sigma,\tau\rangle}$.

Recall that we have $r=a-t+r_1$ for some $r_1\geq 1$. Therefore,
$\zeta^{k^i-k^{i-1}}=\xi^{p^{r_1}}$ for $0\leq i\leq p^t-1$. Let
$\zeta_3\in K$ be a primitive $p^{b+t-\beta}$-th root of unity such
that $\zeta=\zeta_3^{p^{b+t-\beta-a}}$. Then
$\zeta_2=\zeta^{p^{a-b+\beta}}=\zeta_3^{p^t}$.

For $1\leq i\leq p^t-1$ define $w_i=W_i/\zeta_3$ and define
$u_i=U_i/w_i^{p^{r_1}}$. It follows that {\allowdisplaybreaks
\begin{align*}
\sigma\ :\ &u_i\mapsto u_i,~ V_j\mapsto V_j,~ w_i\mapsto\xi w_i,\\
\tau\ :\ &u_1\mapsto u_2\mapsto\cdots\mapsto u_{p^t-1}\mapsto (u_1\cdots u_{p^t-1})^{-1},\\
&V_0\mapsto V_1\mapsto\cdots\mapsto V_{p^t-1}\mapsto\zeta_2^{p^{t-\beta}} V_0,\\
&w_1\mapsto w_2\mapsto\cdots\mapsto w_{p^t-1}\mapsto (w_1\cdots
w_{p^t-1})^{-1},
\end{align*}}
for $1\leq i\leq p^t-1,0\leq j\leq p^t-1$.

Write $L=K(V_j,w_i:1\leq i\leq p^t-1,0\leq j\leq p^t-1)$ and
consider $L(u_i:1\leq i\leq p^t-1)^{\langle\sigma,\tau\rangle}$.
Note that the group $\langle\sigma,\tau\rangle$ acts on the field
$L(u_i)$ as
$\langle\sigma,\tau\rangle/\langle\sigma^{p^t},\tau^{p^b}\rangle$
and is faithful on $L$. According to Lemma \ref{l2.7}, we can
linearize the action of $\tau$ on $u_1,\dots,u_{p^t-1}$. Thus we may
apply Theorem \ref{t2.1} to $L(u_i:1\leq i\leq
p^t-1)^{\langle\sigma,\tau\rangle}$. It remains to show that
$L^{\langle\sigma,\tau\rangle}$ is rational over $K$.

Define $s_1=w_1^{p^t},s_j=w_j/w_{j-1}$ for $2\leq j\leq p^t-1$. We
have $L^{\langle\sigma\rangle}=K(V_i,s_j:0\leq i\leq p^t-1,1\leq
j\leq p^t-1)$. The action of $\tau$ is {\allowdisplaybreaks
\begin{align*}
\tag{7.2} \tau\ :\ &V_0\mapsto V_1\mapsto\cdots\mapsto V_{p^t-1}\mapsto\zeta_2^{p^{t-\beta}} V_0,\\
&s_1\mapsto s_2^{p^t}s_1,~ s_2\mapsto s_3\mapsto\cdots\mapsto
s_{p^t-1}\mapsto 1/(s_1s_2^{p^t-1}\cdots s_{p^t-1}^2).
\end{align*}}
Define $t_1=s_2,t_i=\tau^{i-1}\cdot s_2$ for $2\leq i\leq p^t-1$.
Then $K(t_i,V_j)=K(s_i,V_j)$ and
\begin{align*}
\tau\ :\ &t_1\mapsto t_2\mapsto\cdots\mapsto t_{p^t-1}\mapsto
(t_1t_2\cdots t_{p^t-1})^{-1}.
\end{align*}
The action of $\tau$ on $K(t_i:1\leq i\leq p^t-1)$ can be linearized
according to Lemma \ref{l2.7}. Since $\tau$ acts faithfully on
$K(V_j)$, we can apply Theorem \ref{t2.1}. It remains to show that
$K(V_j)^{\langle\tau\rangle}$ is rational over $K$.

Note that $\tau^{p^t}\cdot V_i=\zeta_2^{p^{t-\beta}}V_i$ for all
$i$, where $\zeta_2^{p^{t-\beta}}=\zeta^{p^{a+t-b}}$ is a primitive
$p^{b-t}$-th root of unity. Define
$v_0=V_0^{p^{b-t}},v_i=V_i/V_{i-1}$ for $1\leq i\leq p^t-1$. Then
$K(v_i)=K(V_i)^{\langle\tau^{p^t}\rangle}$ and we have
\begin{align*}
\tau\ :\ &v_0\mapsto v_0v_1^{p^{b-t}},~ v_1\mapsto
v_2\mapsto\cdots\mapsto v_{p^t-1}\mapsto\zeta_2^{p^{t-\beta}}
(v_1v_2\cdots v_{p^t-1})^{-1}.
\end{align*}
Put $\zeta_3=\zeta^{p^{a-b}}$, a primitive $p^b$-th root of unity.
Whence $\zeta_2^{p^{t-\beta}}=\zeta_3^{p^t}$. Define
$r_i=v_i/\zeta_3$ for $1\leq i\leq p^t-1$. It follows that
\begin{align*}
\tau\ :\ &r_1\mapsto r_2\mapsto\cdots\mapsto r_{p^t-1}\mapsto\
(r_1r_2\cdots r_{p^t-1})^{-1}.
\end{align*}
Applying Lemma \ref{l2.7} once again, we are done.

{\bf Subcase 2.b.} $b-\beta< t$. Define
$\rho_1=\tau^{p^t}\rho^{-p^{\beta-b+t}}$. Then $H$ is isomorphic to
the direct product
$\langle\sigma\rangle\times\langle\rho\rangle\times\langle\rho_1\rangle$,
where $\sigma^{p^a}=\rho^{p^t}=\rho_1^{p^{b-t}}=1$. We have
$\tau^{p^t}=\rho_1\rho^{p^{\beta-b+t}}$.

Define $X_1,X_2\in V^*$ by
\begin{equation*}
X_1=\sum_{i,j}x(\rho^i\rho_1^j),~ X_2=\sum_{i,j}x(\sigma^i\rho^j),~
X_3=\sum_{i,j}x(\sigma^i\rho_1^j).
\end{equation*}
Note that $\sigma\cdot X_2=X_2,\sigma\cdot X_3=X_3,\rho\cdot
X_1=X_1,\rho\cdot X_2=X_2,\rho_1\cdot X_1=X_1$ and $\rho_1\cdot
X_3=X_3$.

Let $\zeta=\zeta_{p^a}\in K$ be a primitive $p^a$-th root of unity.
Define $\zeta_1=\zeta^{p^{a-b+t}}$, a primitive $p^{b-t}$-th root of
unity; and $\xi=\zeta^{p^{a-t}}$, a primitive $p^t$-th root of
unity.

Define $Y_1,Y_2,Y_3\in V^*$ by
\begin{equation*}
Y_1=\sum_{i=0}^{p^a-1}\zeta^{-i}\sigma^i\cdot X_1,~
Y_2=\sum_{i=0}^{p^t-1}\xi^{-i}\rho^i\cdot X_3,~
Y_3=\sum_{i=0}^{p^{b-t}-1}\zeta_1^{-i}\rho_1^i\cdot X_2.
\end{equation*}

It follows that {\allowdisplaybreaks \begin{align*}
\sigma\ :\ &Y_1\mapsto\zeta Y_1,~ Y_2\mapsto Y_2,~ Y_3\mapsto Y_3,\\
\rho\ :\ &Y_1\mapsto Y_1,~ Y_2\mapsto\xi Y_2,~ Y_3\mapsto Y_3,\\
\rho_1\ :\ &Y_1\mapsto Y_1,~ Y_2\mapsto Y_2,~ Y_3\mapsto\zeta_1 Y_3.
\end{align*}}
Thus $K\cdot Y_1+K\cdot Y_2+K\cdot Y_3$ is a representation space of
the subgroup $H$.

Define $x_i=\tau^i\cdot Y_1,y_i=\tau^i\cdot Y_2,z_i=\tau^i\cdot Y_3$
for $0\leq i\leq p^t-1$. We have now {\allowdisplaybreaks
\begin{align*}
\sigma\ :\ &x_i\mapsto \zeta^{k^i} x_i,~ y_i\mapsto\xi^iy_i,~ z_i\mapsto z_i,\\
\tau\ :\ &x_0\mapsto x_1\mapsto\cdots\mapsto x_{p^t-1}\mapsto x_0,\\
&y_0\mapsto y_1\mapsto\cdots\mapsto y_{p^t-1}\mapsto\xi^{p^{\beta-b+t}} y_0,\\
&z_0\mapsto z_1\mapsto\cdots\mapsto z_{p^t-1}\mapsto\zeta_1 z_0,\\
\rho\ :\ &x_i\mapsto x_i,~ y_i\mapsto\xi y_i,~ z_i\mapsto z_i,
\end{align*}}
for $0\leq i\leq p^t-1$. We find that $Y=(\bigoplus_{0\leq i\leq
p^t-1}K\cdot x_i)\oplus(\bigoplus_{0\leq i\leq p^t-1}K\cdot
y_i)\oplus(\bigoplus_{0\leq i\leq p^t-1}K\cdot z_i)$ is a faithful
$\widetilde G$-subspace of $V^*$. Thus, by Theorem \ref{t2.1}, it
suffices to show that $K(x_i,y_i,z_i:0\leq i\leq p^t-1)^{\widetilde
G}$ is rational over $K$.

For $1\leq i\leq p^t-1$, define $U_i=x_i/x_{i-1}$ and
$V_i=y_i/y_{i-1}$. Thus $K(x_i,y_i,z_i:0\leq i\leq
p^t-1)=K(x_0,y_0,U_i,V_i,z_j:1\leq i\leq p^t-1,0\leq j\leq p^t-1)$
and for every $g\in\widetilde G$
\begin{equation*}
g\cdot x_0\in K(U_i,V_i,z_j)\cdot x_0,~ g\cdot y_0\in
K(U_i,V_i,z_j)\cdot y_0,
\end{equation*}
while the subfield $K(U_i,V_i,z_j:1\leq i\leq p^t-1,0\leq j\leq
p^t-1)$ is invariant by the action of $\widetilde G$, i.e.,
{\allowdisplaybreaks
\begin{align*}
\sigma\ :\ &U_i\mapsto \zeta^{k^i-k^{i-1}} U_i,~ V_i\mapsto\xi V_i,~ z_j\mapsto z_j,\\
\tau\ :\ &U_1\mapsto U_2\mapsto\cdots\mapsto U_{p^t-1}\mapsto (U_1\cdots U_{p^t-1})^{-1},\\
\tag{7.3} &V_1\mapsto V_2\mapsto\cdots\mapsto V_{p^t-1}\mapsto\xi^{p^{\beta-b+t}} (V_1\cdots V_{p^t-1})^{-1},\\
&z_0\mapsto z_1\mapsto\cdots\mapsto z_{p^t-1}\mapsto\zeta_1 z_0,\\
\rho\ :\ &U_i\mapsto U_i,~ V_i\mapsto V_i,~ z_j\mapsto z_j,
\end{align*}}
for $1\leq i\leq p^t-1,0\leq j\leq p^t-1$. From Theorem \ref{t2.2}
it follows that if $K(U_i,V_i,z_j:1\leq i\leq p^t-1,0\leq j\leq
p^t-1)^{\widetilde G}$ is rational over $K$, so is
$K(x_i,y_i,z_i:0\leq i\leq p^t-1)^{\widetilde G}$ over $K$.

Compare the actions of $\sigma,\tau,\rho$ in Formula (7.3) with
those in Formula (7.1). They look almost the same. Let $\zeta_2$ be
a primitive $p^{b+t-\beta}$-th root of unity such that
$\zeta=\zeta_2^{p^{b+t-\beta-a}}$. This follows from the assumption
that $K$ contains a primitive root of unity of degree the exponent
of $\widetilde G$. Then $\zeta_2^{p^t}=\xi^{p^{\beta-b+t}}$. Define
$v_i=V_i/\zeta_2$ and use the same method in Subcase 2.a.

{\bf Case 3.} $\widetilde G=G_3$. We have the relations
$(\sigma^{p^{a-\alpha}}\rho^{-1})^{p^\alpha}=1$ and
$\tau^{-1}\sigma^{p^{a-\alpha}}\tau=\sigma^{p^{a-\alpha}}\rho^{p^{a-\alpha}}$.
Clearly, $\widetilde
G'\cap\langle\sigma^{p^{a-\alpha}}\rho^{-1}\rangle=\{1\}$. If
$a-\alpha\geq t$ then $\sigma^{p^{a-\alpha}}\in Z(\widetilde G)$.
Theorem \ref{t2.8} then implies that we can reduce the rationality
problem of $K(\widetilde G)$ to $K(\widetilde
G/\langle\sigma^{p^{a-\alpha}}\rho^{-1}\rangle)$ over $K$, where
$\widetilde G/\langle\sigma^{p^{a-\alpha}}\rho^{-1}\rangle$ clearly
is a metacyclic $p$-group. Therefore, we will assume henceforth that
$a-\alpha<t$.

We are going to show now that $\tau^{p^t}\in Z(\widetilde G)$.
Indeed, we have
\begin{equation*}
\tau^{-p^t}\sigma\tau^{p^t}=\sigma^{(1+p^r)^{p^t}}=\sigma\cdot\sigma^{Ap^{r+t}}=\sigma\cdot\rho^{Ap^{\alpha+r+t-a}}=\sigma,
\end{equation*}
where $A$ is some integer and $\alpha+r+t-a>t$ since $a<\alpha+t\leq
\alpha+r$.

In this way, we obtain that the subgroup
$H=\langle\sigma,\rho,\tau^{p^t}\rangle$ is abelian and has an order
$p^{a+b}$. Put $\rho_1=\sigma^{p^{a-\alpha}}\rho^{-1}$ and
$\rho_2=\tau^{p^t}$. Then
$H\cong\langle\sigma\rangle\times\langle\rho_1\rangle\times\langle\rho_2\rangle$,
where
$\sigma^{p^{a+t-\alpha}}=\rho_1^{p^\alpha}=\rho_2^{p^{b-t}}=1$. We
have also $\rho=\sigma^{p^{a-\alpha}}\rho_1^{-1}$.

Define $X_1,X_2,X_3\in V^*$ by
\begin{equation*}
X_1=\sum_{i,j}x(\rho_1^i\rho_2^j),~
X_2=\sum_{i,j}x(\sigma^i\rho_1^j),~
X_3=\sum_{i,j}x(\sigma^i\rho_2^j).
\end{equation*}
Note that $\sigma\cdot X_2=X_2,\sigma\cdot X_3=X_3,\rho_1\cdot
X_1=X_1,\rho_1\cdot X_2=X_2,\rho_2\cdot X_1=X_1$ and $\rho_2\cdot
X_3=X_3$.

Let $\zeta_1\in K$ be a primitive $p^{a+t-\alpha}$-th root of unity.
Define $\zeta_3=\zeta_1^{p^{a+2t-\alpha-b}}$, a primitive
$p^{b-t}$-th root of unity; and $\xi=\zeta_1^{p^{a-\alpha}}$, a
primitive $p^t$-th root of unity. Since $a-\alpha<t\leq r$, we can
write $r=a-\alpha+r_1$ for some $r_1>0$. Define
$\zeta_2=\xi^{(1+p^{r_1})p^{t-\alpha}}$, a primitive $p^\alpha$-th
root of unity.

Define $Y_1,Y_2,Y_3\in V^*$ by
\begin{equation*}
Y_1=\sum_{i=0}^{p^{a+t-\alpha}-1}\zeta_1^{-i}\sigma^i\cdot X_1,~
Y_2=\sum_{i=0}^{p^\alpha-1}\zeta_2^{-i}\rho_1^i\cdot X_3,~
Y_3=\sum_{i=0}^{p^{b-t}-1}\zeta_3^{-i}\rho_2^i\cdot X_2.
\end{equation*}

It follows that {\allowdisplaybreaks \begin{align*}
\sigma\ :\ &Y_1\mapsto\zeta_1 Y_1,~ Y_2\mapsto Y_2,~ Y_3\mapsto Y_3,\\
\rho_1\ :\ &Y_1\mapsto Y_1,~ Y_2\mapsto\zeta_2 Y_2,~ Y_3\mapsto Y_3,\\
\rho_2\ :\ &Y_1\mapsto Y_1,~ Y_2\mapsto Y_2,~ Y_3\mapsto\zeta_3 Y_3,\\
\rho\ :\ &Y_1\mapsto\xi Y_1,~ Y_2\mapsto\zeta_2^{-1} Y_2,~
Y_3\mapsto Y_3.
\end{align*}}
Thus $K\cdot Y_1+K\cdot Y_2+K\cdot Y_3$ is a representation space of
the subgroup $H$.

Define $x_i=\tau^i\cdot Y_1,y_i=\tau^i\cdot Y_2,z_i=\tau^i\cdot Y_3$
for $0\leq i\leq p^t-1$. We have now {\allowdisplaybreaks
\begin{align*}
\sigma\ :\ &x_i\mapsto \zeta_1^{k^i}\xi^i x_i,~ y_i\mapsto\zeta_2^{-i}y_i,~ z_i\mapsto z_i,\\
\tau\ :\ &x_0\mapsto x_1\mapsto\cdots\mapsto x_{p^t-1}\mapsto x_0,\\
&y_0\mapsto y_1\mapsto\cdots\mapsto y_{p^t-1}\mapsto y_0,\\
&z_0\mapsto z_1\mapsto\cdots\mapsto z_{p^t-1}\mapsto\zeta_3 z_0,\\
\rho\ :\ &x_i\mapsto\xi x_i,~ y_i\mapsto\zeta_2^{-1} y_i,~
z_i\mapsto z_i.
\end{align*}}
for $0\leq i\leq p^t-1$. We find that $Y=(\bigoplus_{0\leq i\leq
p^t-1}K\cdot x_i)\oplus(\bigoplus_{0\leq i\leq p^t-1}K\cdot
y_i)\oplus(\bigoplus_{0\leq i\leq p^t-1}K\cdot z_i)$ is a faithful
$\widetilde G$-subspace of $V^*$. Thus, by Theorem \ref{t2.1}, it
suffices to show that $K(x_i,y_i,z_i:0\leq i\leq p^t-1)^{\widetilde
G}$ is rational over $K$.

For $1\leq i\leq p^t-1$, define $U_i=x_i/x_{i-1}$ and
$V_i=y_i/y_{i-1}$. Thus $K(x_i,y_i,z_i:0\leq i\leq
p^t-1)=K(x_0,y_0,U_i,V_i,z_j:1\leq i\leq p^t-1,0\leq j\leq p^t-1)$
and for every $g\in\widetilde G$
\begin{equation*}
g\cdot x_0\in K(U_i,V_i,z_j)\cdot x_0,~ g\cdot y_0\in
K(U_i,V_i,z_j)\cdot y_0,
\end{equation*}
while the subfield $K(U_i,V_i,z_j:1\leq i\leq p^t-1,0\leq j\leq
p^t-1)$ is invariant by the action of $\widetilde G$, i.e.,
{\allowdisplaybreaks
\begin{align*}
\sigma\ :\ &U_i\mapsto \zeta_1^{k^i-k^{i-1}}\xi U_i,~ V_i\mapsto\zeta_2^{-1} V_i,~ z_j\mapsto z_j,\\
\tau\ :\ &U_1\mapsto U_2\mapsto\cdots\mapsto U_{p^t-1}\mapsto (U_1\cdots U_{p^t-1})^{-1},\\
\tag{7.4} &V_1\mapsto V_2\mapsto\cdots\mapsto V_{p^t-1}\mapsto (V_1\cdots V_{p^t-1})^{-1},\\
&z_0\mapsto z_1\mapsto\cdots\mapsto z_{p^t-1}\mapsto\zeta_3 z_0,\\
\rho\ :\ &U_i\mapsto U_i,~ V_i\mapsto V_i,~ z_j\mapsto z_j.
\end{align*}}
for $1\leq i\leq p^t-1$ and $0\leq j\leq p^t-1$. From Theorem
\ref{t2.2} it follows that if $K(U_i,V_i,z_j:1\leq i\leq p^t-1,0\leq
j\leq p^t-1)^{\widetilde G}$ is rational over $K$, so is
$K(x_i,y_i,z_i:0\leq i\leq p^t-1)^{\widetilde G}$ over $K$.

Since $\rho$ acts trivially on $K(U_i,V_i,z_j)$, we find that
$K(U_i,V_i,z_j)^{\widetilde
G}=K(U_i,V_i,z_j)^{\langle\sigma,\tau\rangle}$.

Recall that $r_1=r-a+\alpha>0$. Therefore,
$\zeta_1^{k-1}=\zeta_1^{p^{a-\alpha+r_1}}=\xi^{p^{r_1}}$ and also
$\zeta_1^{k^{i-1}(k-1)}=\xi^{p^{r_1}}$ for all $i$.

Define $v_i=U_i^{p^{t-\alpha}}V_i$. Since
$\zeta_2=\xi^{(1+p^{r_1})p^{t-\alpha}}$, we have
{\allowdisplaybreaks
\begin{align*}
\sigma\ :\ &U_i\mapsto \xi^{1+p^{r_1}} U_i,~ v_i\mapsto v_i,~ z_j\mapsto z_j,\\
\tau\ :\ &U_1\mapsto U_2\mapsto\cdots\mapsto U_{p^t-1}\mapsto (U_1\cdots U_{p^t-1})^{-1},\\
&v_1\mapsto v_2\mapsto\cdots\mapsto v_{p^t-1}\mapsto (v_1\cdots v_{p^t-1})^{-1},\\
&z_0\mapsto z_1\mapsto\cdots\mapsto z_{p^t-1}\mapsto\zeta_3 z_0.
\end{align*}}
for $1\leq i\leq p^t-1$ and $0\leq j\leq p^t-1$.

Write $L=K(U_i,z_j:1\leq i\leq p^t-1,0\leq j\leq p^t-1)$ and
consider $L(v_i:1\leq i\leq p^t-1)^{\langle\sigma,\tau\rangle}$.
Note that the group $\langle\sigma,\tau\rangle$ acts on the field
$L(v_i)$ as $\langle\sigma,\tau\rangle/\langle\sigma^{p^t}\rangle$
and is faithful on $L$. According to Lemma \ref{l2.7}, we can
linearize the action of $\tau$ on $v_1,\dots,v_{p^t-1}$. Thus we may
apply Theorem \ref{t2.1} to $L(v_i:1\leq i\leq
p^t-1)^{\langle\sigma,\tau\rangle}$. It remains to show that
$L^{\langle\sigma,\tau\rangle}$ is rational over $K$.

Define $u_1=U_1^{p^t},u_i=U_i/U_{i-1}$ for $2\leq i\leq p^t-1$. Then
$K(u_i,z_j,1\leq i\leq p^t-1,0\leq j\leq
p^t-1)=L^{\langle\sigma\rangle}$ and the action of $\tau$ is
{\allowdisplaybreaks
\begin{align*}
\tag{7.5} \tau\ :\ &u_1\mapsto u_1u_2^{p^t}, u_2\mapsto u_3\mapsto\cdots\mapsto (u_1u_2^{p^t-1}u_3^{p^t-2}\cdots u_{p^t-1}^2)^{-1},\\
&z_0\mapsto z_1\mapsto\cdots\mapsto z_{p^t-1}\mapsto\zeta_3 z_0.
\end{align*}}
Compare the action of $\tau$ in Formula (7.5) with that in Formula
(7.2). Use the same method in Subcase 2.a.

{\bf Case 4.} $\widetilde G=G_4$. As in Case 3, we may assume that
$a<\alpha+t$. Therefore, $\tau^{p^t}$ is in the center of
$\widetilde G$ and the group
$H=\langle\sigma,\rho,\tau^{p^t}\rangle$ is abelian of order
$p^{a+b}$. There are two possibilities for the decomposition of $H$
as a direct product of cyclic groups.

{\bf Subcase 4.a.} $b-\beta\leq t$. Define
$\rho_1=\sigma^{p^{a-\alpha}}\rho^{-1},\rho_2=\tau^{p^t}\rho^{-p^{\beta+t-b}}$.
Then $H$ is isomorphic to the direct product
$\langle\sigma\rangle\times\langle\rho_1\rangle\times\langle\rho_2\rangle$,
where
$\sigma^{p^{a+t-\alpha}}=\rho_1^{p^\alpha}=\rho_2^{p^{b-t}}=1$. We
have $\rho=\rho_1^{-1}\sigma^{p^{a-\alpha}}$ and
$\tau^{p^t}=\rho_2\rho^{p^{\beta+t-b}}$.

We can adopt exactly the same definitions from Case 3, starting with
$X_1,X_2.X_3\in V^*$ till the place where we define $x_i=\tau^i\cdot
Y_1,y_i=\tau^i\cdot Y_2,z_i=\tau^i\cdot Y_3$ for $0\leq i\leq
p^t-1$. Here only the action of $\tau$ is changed a little:
{\allowdisplaybreaks
\begin{align*}
\sigma\ :\ &x_i\mapsto \zeta_1^{k^i}\xi^i x_i,~ y_i\mapsto\zeta_2^{-i}y_i,~ z_i\mapsto z_i,\\
\tau\ :\ &x_0\mapsto x_1\mapsto\cdots\mapsto x_{p^t-1}\mapsto\xi^{p^{\beta+t-b}} x_0,\\
&y_0\mapsto y_1\mapsto\cdots\mapsto y_{p^t-1}\mapsto\zeta_2^{-p^{\beta+t-b}} y_0,\\
&z_0\mapsto z_1\mapsto\cdots\mapsto z_{p^t-1}\mapsto\zeta_3 z_0,\\
\rho\ :\ &x_i\mapsto\xi x_i,~ y_i\mapsto\zeta_2^{-1} y_i,~
z_i\mapsto z_i,
\end{align*}}
for $0\leq i\leq p^t-1$.

Let $\zeta_4\in K$ be a primitive $p^{b+t-\beta}$-th root of unity
such that $\zeta_4=\zeta_1^{p^{a-\alpha-b+\beta}}$. Then
$\zeta_4^{p^t}=\xi^{p^{\beta-b+t}}$. Similarly, we can define
$p^{\alpha+b-\beta}$-th root of unity $\zeta_5$ such that
$\zeta_5=(\xi^{1+p^{r_1}})^{p^{t-(b-\beta+\alpha)}}$. Then
$\zeta_2^{p^{\beta-b+t}}=\zeta_5^{p^t}$. Note that
$\alpha+b-\beta<t+b-\beta$, so $\zeta_5$ is contained in $K$.

Define $U_i=x_i/\zeta_4x_{i-1}$ and $V_i=\zeta_5y_i/y_{i-1}$ for
$1\leq i\leq p^t-1$. Thus we obtain exactly the same actions in Case
3 given by Formula (7.4).

{\bf Subcase 4.b.} $b-\beta> t$. We have
$(\tau^{p^{b-\beta}}\rho^{-1})^{p^\beta}=1$. Then $\widetilde G'\cap
\langle\tau^{p^{b-\beta}}\rho^{-1}\rangle=\{1\}$, so we can apply
Theorem \ref{t2.8} reducing the rationality problem of $K(\widetilde
G)$ to the rationality problem of $K(\widetilde G_1)$ over $K$,
where $\widetilde G_1\cong\widetilde
G/\langle\tau^{p^{b-\beta}}\rho^{-1}\rangle$. The group $\widetilde
G_1$ is generated by elements $\sigma,\tau$ and $\rho$ such that
$\sigma^{p^a}=\rho^{p^\alpha},\tau^{p^{b-\beta}}=\rho,\rho^{p^t}=1,[\sigma,\tau]=\rho,
\rho - \text{central}$. Then the abelian subgroup
$H=\langle\sigma,\rho,\tau^{p^t}\rangle$ is of order
$p^{a+b-\beta}$. Put $y=b+\alpha-\beta-t$.

{\bf Sub-subcase 4.b.a.} Let $y\leq a$. Define
$\rho_1=\tau^{p^t}\sigma^{-p^{a-y}}$. Then $H$ is isomorphic to the
direct product $\langle\sigma\rangle\times\langle\rho_1\rangle$,
where $\sigma^{p^{a+t-\alpha}}=\rho_1^{p^y}=1$. Note that
$\tau^{p^t}=\rho_1\sigma^{p^{a-y}}$ and
$\rho=\sigma^{p^{a-\alpha}}\rho_1^{p^{y-\alpha}}$.

Define $X_1,X_2\in V^*$ by
\begin{equation*}
X_1=\sum_ix(\rho_1^i),~ X_2=\sum_ix(\sigma^i).
\end{equation*}
Note that $\sigma\cdot X_2=X_2$ and $\rho_1\cdot X_1=X_1$.

Recall that $r_1=r-a+\alpha>0$. Let $\zeta_1\in K$ be a primitive
$p^{a+t-\alpha}$-th root of unity, and put
$\xi=\zeta_1^{p^{a-\alpha}}$, a primitive $p^t$-th root of unity.
Then $\zeta_1^{k-1}=\zeta_1^{p^{a-\alpha+r_1}}=\xi^{p^{r_1}}$ and
also $\zeta_1^{k^{i-1}(k-1)}=\xi^{p^{r_1}}$ for all $i$. Define
$\zeta_2=\zeta_1^{(1+p^{r_1})p^{a+t-\alpha-y}}$, a primitive
$p^y$-th root of unity.

Define $Y_1,Y_2\in V^*$ by
\begin{equation*}
Y_1=\sum_{i=0}^{p^{a+t-\alpha}-1}\zeta_1^{-i}\sigma^i\cdot X_1,~
Y_2=\sum_{i=0}^{p^y-1}\zeta_2^{-i}\rho_1^i\cdot X_2.
\end{equation*}

It follows that {\allowdisplaybreaks \begin{align*}
\sigma\ :\ &Y_1\mapsto\zeta_1 Y_1,~ Y_2\mapsto Y_2,\\
\rho_1\ :\ &Y_1\mapsto Y_1,~ Y_2\mapsto\zeta_2 Y_2,\\
\rho\ :\ &Y_1\mapsto\xi Y_1,~ Y_2\mapsto\zeta_2^{p^{y-\alpha}} Y_2.
\end{align*}}
Thus $K\cdot Y_1+K\cdot Y_2$ is a representation space of the
subgroup $H$.

Define $x_i=\tau^i\cdot Y_1,y_i=\tau^i\cdot Y_2$ for $0\leq i\leq
p^t-1$. We have now {\allowdisplaybreaks
\begin{align*}
\sigma\ :\ &x_i\mapsto \zeta_1^{k^i}\xi^i x_i,~ y_i\mapsto(\zeta_2^{p^{y-\alpha}})^iy_i,\\
\tau\ :\ &x_0\mapsto x_1\mapsto\cdots\mapsto x_{p^t-1}\mapsto\zeta_1^{p^{a-y}} x_0,\\
&y_0\mapsto y_1\mapsto\cdots\mapsto y_{p^t-1}\mapsto\zeta_2 y_0,\\
\rho\ :\ &x_i\mapsto\xi x_i,~ y_i\mapsto\zeta_2^{p^{y-\alpha}} y_i,
\end{align*}}
for $0\leq i\leq p^t-1$.

For $1\leq i\leq p^t-1$ define $U_i=x_i/x_{i-1},V_i=y_i/y_{i-1}$. It
follows that {\allowdisplaybreaks
\begin{align*}
\sigma\ :\ &U_i\mapsto \xi^{1+p^{r_1}}U_i,~ V_i\mapsto\zeta_2^{p^{y-\alpha}} V_i,\\
\tau\ :\ &U_1\mapsto U_2\mapsto\cdots\mapsto U_{p^t-1}\mapsto\zeta_1^{p^{a-y}} (U_1\cdots U_{p^t-1})^{-1},\\
&V_1\mapsto V_2\mapsto\cdots\mapsto V_{p^t-1}\mapsto\zeta_2
(V_1\cdots V_{p^t-1})^{-1},\\
\rho\ :\ &U_i\mapsto U_i,~ V_i\mapsto V_i,
\end{align*}}
for $1\leq i\leq p^t-1$.

Since $\rho$ acts trivially on $K(U_i,V_i:1\leq i\leq p^t-1)$, we
find that $K(U_i,V_i)^{\widetilde
G}=K(U_i,V_i)^{\langle\sigma,\tau\rangle}$.

Note that $\zeta_2^{p^{y-\alpha}}=\xi^{(1+p^{r_1})p^{t-\alpha}}$.
For $1\leq i\leq p^t-1$ define $v_i=U_i^{-p^{t-\alpha}}V_i$. We have
{\allowdisplaybreaks
\begin{align*}
\sigma\ :\ &U_i\mapsto \xi^{1+p^{r_1}}U_i,~ v_i\mapsto v_i,\\
\tau\ :\ &U_1\mapsto U_2\mapsto\cdots\mapsto U_{p^t-1}\mapsto\zeta_1^{p^{a-y}} (U_1\cdots U_{p^t-1})^{-1},\\
&v_1\mapsto v_2\mapsto\cdots\mapsto
v_{p^t-1}\mapsto\zeta_1^{p^{r+t-y}} (v_1\cdots v_{p^t-1})^{-1},
\end{align*}}
for $1\leq i\leq p^t-1$.

Let $\zeta_3\in K$ be a primitive $p^{b+t-\beta}$-th root of unity
such that $\zeta_1=\zeta_3^{p^{b+t-\beta-(a+t-\alpha)}}$. Then
$\zeta_3^{p^t}=\zeta_1^{p^{a-y}}$. Define $u_1=U_1^{p^t}$ and for
$2\leq i\leq p^t-1$ define $u_i=U_i/(\zeta_3U_{i-1})$.

From the inequality $a<\alpha+t$ we get $a-r-t<0$. Whence
$y-r_1+t=b-\beta+t+a-r-t<b-\beta+t$. Let $\zeta_4\in K$ be a
primitive $p^{y-r_1+t}$-th root of unity such that
$\zeta_4=\zeta_3^{p^{b+t-\beta-(y-r_1+t)}}$. Then
$\zeta_4^{p^t}=\zeta_1^{p^{r+t-y}}$.

For $1\leq i\leq p^t-1$ define $w_i=v_i/\zeta_4$. The actions of
$\sigma$ and $\tau$ are then given by {\allowdisplaybreaks
\begin{align*}
\sigma\ :\ &u_i\mapsto u_i,~ w_i\mapsto w_i,\\
\tau\ :\ &u_1\mapsto u_1u_2^{p^t},~ u_2\mapsto u_3\mapsto\cdots\mapsto (u_1u_2^{p^t-1}u_3^{p^t-2}\cdots u_{p^t-1}^2)^{-1},\\
&w_1\mapsto w_2\mapsto\cdots\mapsto w_{p^t-1}\mapsto (w_1\cdots
w_{p^t-1})^{-1},
\end{align*}}
for $1\leq i\leq p^t-1$. Since $\sigma$ acts trivially on
$K(u_i,w_i:1\leq i\leq p^t-1)$, we find that
$K(u_i,w_i)^{\langle\sigma,\tau\rangle}=K(u_i,w_i)^{\langle\tau\rangle}$.
As we did before, we can easily linearize the action of $\tau$
applying Lemma \ref{l2.7}.

{\bf Sub-subcase 4.b.b.} Let $y>a$. Put $u=y-a>0$. Define
$\rho_1=\sigma\tau^{-p^{t+u}}$ and $\rho_2=\tau^{p^t}$. Therefore,
$H\cong\langle\rho_1\rangle\times\langle\rho_2\rangle$, where
$\rho_1^{p^a}=\rho_2^{p^{b-\beta}}=1,\sigma=\rho_1\rho_2^{p^u}$ and
$\rho=\rho_2^{p^{b-\beta-t}}$.

Define $X_1,X_2\in V^*$ by
\begin{equation*}
X_1=\sum_ix(\rho_1^i),~ X_2=\sum_ix(\rho_2^i).
\end{equation*}
Note that $\rho_1\cdot X_1=X_1$ and $\rho_2\cdot X_2=X_2$.

Let $\zeta_2\in K$ be a primitive $p^{b-\beta}$-th root of unity.
Put $\zeta_1=\zeta_2^{p^u}$, a primitive $p^{a+t-\alpha}$-th root of
unity; $\xi=\zeta_1^{p^{a-\alpha}}$, a primitive $p^t$-th root of
unity; $\zeta_3=\zeta_1^{p^{t-\alpha}}$, a primitive $p^a$-th root
of unity.

Define $Y_1,Y_2\in V^*$ by
\begin{equation*}
Y_1=\sum_{i=0}^{p^a-1}\zeta_3^{-i}\rho_1^i\cdot X_2,~
Y_2=\sum_{i=0}^{p^{b-\beta}-1}\zeta_2^{-i}\rho_2^i\cdot X_1.
\end{equation*}

It follows that {\allowdisplaybreaks \begin{align*}
\rho_1\ :\ &Y_1\mapsto\zeta_3 Y_1,~ Y_2\mapsto Y_2,\\
\rho_2\ :\ &Y_1\mapsto Y_1,~ Y_2\mapsto\zeta_2 Y_2,\\
\sigma\ :\ &Y_1\mapsto\zeta_3 Y_1,~ Y_2\mapsto\zeta_1 Y_2,\\
\rho\ :\ &Y_1\mapsto Y_1,~ Y_2\mapsto\xi Y_2.
\end{align*}}
Thus $K\cdot Y_1+K\cdot Y_2$ is a representation space of the
subgroup $H$.

Define $x_i=\tau^i\cdot Y_1,y_i=\tau^i\cdot Y_2$ for $0\leq i\leq
p^t-1$. We have now {\allowdisplaybreaks
\begin{align*}
\sigma\ :\ &x_i\mapsto \zeta_3^{k^i} x_i,~ y_i\mapsto\zeta_1^{k^i}\xi^iy_i,\\
\tau\ :\ &x_0\mapsto x_1\mapsto\cdots\mapsto x_{p^t-1}\mapsto x_0,\\
&y_0\mapsto y_1\mapsto\cdots\mapsto y_{p^t-1}\mapsto\zeta_2 y_0,\\
\rho\ :\ &x_i\mapsto x_i,~ y_i\mapsto\xi y_i,
\end{align*}}
for $0\leq i\leq p^t-1$. Define $U_i=x_i/x_{i-1}$ and
$V_i=y_i/y_{i-1}$ for $1\leq i\leq p^t-1$. It follows that
{\allowdisplaybreaks
\begin{align*}
\sigma\ :\ &U_i\mapsto\zeta_3^{k^i-k^{i-1}} U_i,~ V_i\mapsto\zeta_1^{k^i-k^{i-1}}\xi V_i,\\
\tau\ :\ &U_1\mapsto U_2\mapsto\cdots\mapsto U_{p^t-1}\mapsto (U_1U_2\cdots U_{p^t-1})^{-1},\\
&V_1\mapsto V_2\mapsto\cdots\mapsto V_{p^t-1}\mapsto\zeta_2 (V_1V_2\cdots V_{p^t-1})^{-1},\\
\rho\ :\ &U_i\mapsto U_i,~ V_i\mapsto V_i,
\end{align*}}
for $0\leq i\leq p^t-1$. Since $\rho$ acts trivially on
$K(U_i,V_i:1\leq i\leq p^t-1)$, we find that $K(U_i,V_i)^{\widetilde
G}=K(U_i,V_i)^{\langle\sigma,\tau\rangle}$.

Recall that $r_1=r-a+\alpha>0$ and that
$\zeta_1^{k^{i-1}(k-1)}=\xi^{p^{r_1}}=\xi^{p^{r+\alpha-a}}$ for all
$i$. We also have
$\zeta_3^{k^{i-1}(k-1)}=\xi^{p^{r_1+t-\alpha}}=\xi^{p^{r+t-a}}$ for
all $i$. Put $\eta=\xi^{1+p^{r+\alpha-a}}$. Since
$\langle\eta\rangle=\langle\xi\rangle$, there exists $w\in\mathbb Z$
such that $\xi^{p^{r+t-a}}=\eta^w$.

For $1\leq i\leq p^t-1$ define $u_i=U_iV_i^{-w}$. It follows that
{\allowdisplaybreaks
\begin{align*}
\sigma\ :\ &u_i\mapsto u_i,~ V_i\mapsto\xi^{1+p^{r+\alpha-a}} V_i\\
\tau\ :\ &u_1\mapsto u_2\mapsto\cdots\mapsto u_{p^t-1}\mapsto\zeta_2^{-w} (u_1u_2\cdots u_{p^t-1})^{-1},\\
&V_1\mapsto V_2\mapsto\cdots\mapsto V_{p^t-1}\mapsto\zeta_2
(V_1V_2\cdots V_{p^t-1})^{-1},
\end{align*}}
for $1\leq i\leq p^t-1$. Similar actions of $\sigma$ and $\tau$ have
been considered many times in the previous cases, so we omit the
details of the final stage of the proof.

{\bf Case 5.} $\widetilde G=G_5$. From Proposition \ref{p3.3} it
follows that $\widetilde G$ is cyclic and is generated by
$\sigma^{-2+2^r}\rho$. If $a\geq t+1$ then
$(\sigma^{-2+2^r}\rho)^{2^{a-1}}=1$ and whence $\widetilde
G\cap\langle\rho\rangle=\{1\}$. Theorem \ref{t2.8} then implies that
we can reduce the rationality problem of $K(\widetilde G)$ to
$K(\widetilde G/\langle\rho\rangle)$ over $K$, where $\widetilde
G/\langle\rho\rangle$ clearly is a metacyclic $p$-group.

Now, let $a=t$. Whence $r=t$, i.e., we have the relations
$\sigma^{2^t}=\rho^{2^t}=\tau^{2^b}=1$ and
$\tau^{-1}\sigma\tau=\sigma^{-1}\rho$. If $b\geq t+1$ we can apply
the proof of Case 5 in Theorem \ref{t1.5}, where we require a
primitive $2^{t+1}$-th root of unity in $K$. Therefore, we will
assume that $\tau^{2^t}=1$.

From $\tau^{-2}\sigma\tau^2=\sigma$ it follows that $\tau^2$ is in
the center $Z(\widetilde G)$. Since $\widetilde
G\cap\langle\tau^2\rangle=\{1\}$, we can apply Theorem \ref{t2.8} so
that we reduce the rationality problem of $K(\widetilde G)$ to
$K(\widetilde G/\langle\tau^2\rangle)$ over $K$. In this way, we can
assume that $\tau^2=1$. Note that
$\ord(\tau^2)=2^{t-1}<\ord(\sigma)=2^t=\exp(\widetilde
G_1/\widetilde G_1')$ for $\widetilde G_1=\widetilde
G/\langle\tau^2\rangle$.

Let $\xi$ be a primitive $p^t$-th root of unity. Define $X_1,X_2\in
V^*$ by
\begin{equation*}
X_1=\sum_{i=0}^{2^t-1}x(\sigma^i),~ X_2=\sum_{i=0}^{2^t-1}x(\rho^i).
\end{equation*}
Note that $\sigma\cdot X_1=X_1$ and $\rho\cdot X_2=X_2$.

Define $Y_1,Y_2\in V^*$ by
\begin{equation*}
Y_1=\sum_{i=0}^{2^t-1}\xi^{-i}\rho^i\cdot X_1,~
Y_2=\sum_{i=0}^{2^t-1}\xi^{-i}\sigma^i\cdot X_2.
\end{equation*}

It follows that {\allowdisplaybreaks \begin{align*}
\sigma\ :\ &Y_1\mapsto Y_1,Y_2\mapsto\xi Y_2,\\
\rho\ :\ &Y_1\mapsto\xi Y_1,Y_2\mapsto Y_2.
\end{align*}}
Thus $K\cdot Y_1+K\cdot Y_2$ is a representation space of the
subgroup $\langle\sigma,\rho\rangle$.

Define $x_i=\tau^i\cdot Y_1,y_i=\tau^i\cdot Y_2$ for $0\leq i\leq
1$. We have now {\allowdisplaybreaks \begin{align*}
\sigma\ :\ &x_i\mapsto \xi^{\delta_i} x_i,~ y_i\mapsto\xi^{\gamma_i}y_i\\
\tau\ :\ &x_0\mapsto x_1\mapsto x_0,\\
&y_0\mapsto y_1\mapsto y_0,\\
\rho\ :\ &x_i\mapsto\xi x_i,~ y_i\mapsto y_i,
\end{align*}}
where $\delta_i=0,\gamma_i=1$ for $i$--even;
$\delta_i=1,\gamma_i=-1$ for $i$--odd; and $0\leq i\leq 1$.

We find that $Y=\bigoplus_{0\leq i\leq 1}K\cdot x_i$ is a faithful
$\widetilde G$-subspace of $V^*$. Thus, by Theorem \ref{t2.1}, it
suffices to show that $K(x_i:0\leq i\leq 1)^{\widetilde G}$ is
rational over $K$. The rationality of $K(x_i:0\leq i\leq
1)^{\widetilde G}$ over $K$ follows from Theorem \ref{t2.9}.

{\bf Case 6.} $\widetilde G=G_6$. Similarly to Case 5, we may assume
that $a=t$. Whence $r=t$, i.e., we have the relations
$\sigma^{2^t}=\rho^{2^t}=1, \tau^{2^b}=\rho^{2^\beta}$ and
$\tau^{-1}\sigma\tau=\sigma^{-1}\rho$. Since in the proof of Case 6
in Theorem \ref{t1.5} we require a primitive $2^{t+1}$-th root of
unity in $K$, and $\ord(\tau)\geq 2^{b+1}\geq 2^{t+1}$, we may apply
the same proof.

{\bf Case 7.} $\widetilde G=G_7$. From the relation
$\tau^{-1}\sigma^{2^a}\tau=\sigma^{-2^a+2^{a+r}}=\sigma^{2^a}$ it
follows that the order of $\sigma$ is $2^{a+1}$. Therefore,
$\alpha=t-1$. We can apply the proof of Case 7 in Theorem
\ref{t1.5}, since there we require only primitive roots of unity of
smaller order than the exponent of $\widetilde G$.

{\bf Case 8.} $\widetilde G=G_8$. The proof is almost the same as in
Case 7. The only difference is the action of $\tau$.

{\bf Case 9.} $\widetilde G=G_9$. Recall that $c\geq r\geq t$. The
proof is exactly the same as in Case 1. Indeed, if $a<r+t$ we have
that $\tau^{p^t}$ is in the center of $\widetilde G$ and $\widetilde
G'\cap\langle\tau^{p^t}\rangle=1$, so we can apply Theorem
\ref{t2.8}. Thus we may assume again that $\tau^{p^t}=1$.

{\bf Case 10.} $\widetilde G=G_{10}$. Similarly to Case 2, we may
assume that $a<t+r$ and $\tau^{p^t}\in Z(\widetilde G)$. The group
$H=\langle\sigma,\rho,\tau^{p^t}\rangle$ is abelian of order
$p^{a+b}$.

{\bf Subcase 10.a.} $b-\beta\geq t$.

{\bf Sub-subcase 10.a.a.} $c+t\geq b$. Define
$\rho_1=\tau^{p^{b-\beta}}\sigma^{-p^{c-\beta}}\rho^{-1}$ and
$\rho_2=\tau^{p^t}\sigma^{-p^{c+t-b}}$. Then $H$ is isomorphic to
the direct product
$\langle\sigma\rangle\times\langle\rho_1\rangle\times\langle\rho_2\rangle$,
where $\sigma^{p^a}=\rho_1^{p^\beta}=\rho_2^{p^{b-\beta}}=1$. We
have $\rho=\rho_1^{-1}\rho_2^{p^{b-\beta-t}}$ and
$\tau^{p^t}=\rho_2\sigma^{p^{c+t-b}}$.

The only difference now with Subcase 2.a is the action of $\tau$ on
$x_i$'s. As we did many times so far, with a proper adjustment of
the variables we can easily repair this 'defect'. A similar
situation will reappear for all remaining cases, subcases and
sub-subcases. For each subcase that follows, it is tiresome but not
difficult to change the variables in such a way that we may directly
apply the proof of the respective subcases of cases 2,3 and 4. As an
illustration, we give some extra explanation in Subcase 11.b. For
the remaining subcases we will write only the decompositions of
$H=\langle\sigma,\rho,\tau^{p^t}\rangle$ as a direct product of
cyclic groups.

{\bf Sub-subcase 10.a.b.} $c+t<b$ and $c+t<a+\beta$. We have
$\tau^{p^{b+t-\beta}}=\sigma^{p^{c+t-\beta}}\ne 1$. If we suppose
that a power of $\sigma^{p^c}\rho^{p^\beta}$ is in $\widetilde G'$,
we get $c=r+\beta$, whence $c+t-\beta=r+t>a$, a contradiction.
Therefore, $\langle\tau^{p^t}\rangle\cap\widetilde G'=\{1\}$ and we
can apply Theorem \ref{t2.8}, reducing this subcase to Case 1.

{\bf Sub-subcase 10.a.c.} $c+t<b$ and $c+t\geq a+\beta$. We have
$\tau^{p^{b+a-c}}=\rho^{p^{\beta+a-c}}$. Define
$\rho_1=\tau^{p^{b-\beta}}\sigma^{-p^{c-\beta}}\rho^{-1}$ and
$\rho_2=\tau^{p^t}$. Then $H$ is isomorphic to the direct product
$\langle\sigma\rangle\times\langle\rho_1\rangle\times\langle\rho_2\rangle$,
where $\sigma^{p^a}=\rho_1^{p^\beta}=\rho_2^{p^{b-\beta}}=1$ and
$\rho=\rho_1^{-1}\rho_2^{p^{b-\beta-t}}\sigma^{-p^{c-\beta}}$.

{\bf Subcase 10.b.} $b-\beta< t$. Then $c\geq\beta>b-t$, i.e.,
$c-b+t>0$. Define
$\rho_1=\tau^{p^t}\sigma^{-p^{c-b+t}}\rho^{-p^{\beta-b+t}}$. Then
$H$ is isomorphic to the direct product
$\langle\sigma\rangle\times\langle\rho\rangle\times\langle\rho_1\rangle$,
where $\sigma^{p^a}=\rho^{p^t}=\rho_1^{p^{b-t}}=1$. We have
$\tau^{p^t}=\rho_1\sigma^{p^{c+t-b}}\rho^{p^{\beta-b+t}}$.

{\bf Case 11.} $\widetilde G=G_{11}$. Similarly to Case 3, we may
assume that $a<t+\alpha$ and $\tau^{p^t}\in Z(\widetilde G)$. The
subgroup $H=\langle\sigma,\rho,\tau^{p^t}\rangle$ is abelian and has
an order $p^{a+b}$.

{\bf Subcase 11.a.} $c-b+t\geq 0$. Put
$\rho_1=\sigma^{p^{a-\alpha}}\rho^{-1}$ and
$\rho_2=\tau^{p^t}\sigma^{-p^{c-b+t}}$. Then
$H\cong\langle\sigma\rangle\times\langle\rho_1\rangle\times\langle\rho_2\rangle$,
where
$\sigma^{p^{a+t-\alpha}}=\rho_1^{p^\alpha}=\rho_2^{p^{b-t}}=1$. We
have also $\rho=\sigma^{p^{a-\alpha}}\rho_1^{-1}$ and
$\tau^{p^t}=\rho_2\sigma^{p^{c-b+t}}$.

{\bf Subcase 11.b.} $c-b+t<0$. Put
$\rho_1=\sigma\tau^{-p^{b-c}},\rho_2=\sigma^{p^{a-\alpha}}\rho^{-1},\rho_3=\tau^{p^t}$.
Then
$H\cong\langle\rho_1\rangle\times\langle\rho_2\rangle\times\langle\rho_3\rangle$,
where $\rho_1^{p^c}=\rho_2^{p^\alpha}=\rho_3^{p^{a+b-c-\alpha}}=1$.
Note that $\sigma=\rho_1\rho_3^{p^{b-c-t}}$ and
$\rho=\rho_1^{p^{a-\alpha}}\rho_2^{-1}\rho_3^{p^{a-\alpha+b-c-t}}$.
In this situation we need to adjust first the variables so that the
actions of $\sigma$ and $\rho$ become the same as in Case 3.

Let $\zeta\in K$ be a primitive $p^{a+b-c-\alpha}$-th root of unity.
Define $\zeta_1=\zeta^{p^{b-c-t}}$, a primitive $p^{a+t-\alpha}$-th
root of unity; $\zeta_2=\zeta_1^{p^{a+t-\alpha-c}}$, a primitive
$p^c$-th root of unity; and $\xi=\zeta_1^{p^{a-\alpha}}$, a
primitive $p^t$-th root of unity. Since $a-\alpha<t\leq r$, we can
write $r=a-\alpha+r_1$ for some $r_1>0$. Define
$\zeta_3=\xi^{(1+p^{r_1})p^{t-\alpha}}$, a primitive $p^\alpha$-th
root of unity.

Define $Y_1,Y_2,Y_3\in V^*$ so that we have the actions
{\allowdisplaybreaks \begin{align*}
\rho_1\ :\ &Y_1\mapsto\zeta_2 Y_1,~ Y_2\mapsto Y_2,~ Y_3\mapsto Y_3,\\
\rho_2\ :\ &Y_1\mapsto Y_1,~ Y_2\mapsto\zeta_3 Y_2,~ Y_3\mapsto Y_3,\\
\rho_3\ :\ &Y_1\mapsto Y_1,~ Y_2\mapsto Y_2,~ Y_3\mapsto\zeta Y_3,\\
\rho\ :\ &Y_1\mapsto\zeta_2^{p^{a-\alpha}} Y_1,~
Y_2\mapsto\zeta_3^{-1} Y_2,~
Y_3\mapsto\xi Y_3,\\
\sigma\ :\ &Y_1\mapsto\zeta_2 Y_1,~ Y_2\mapsto Y_2,~
Y_3\mapsto\zeta_1 Y_3.
\end{align*}}
Define $Z_1=Y_3,Z_2=Y_2,Z_3=Y_3^{p^{a+t-\alpha-c}}Y_1^{-1}$. It is
easily seen now that the actions of $\sigma$ and $\rho$ on
$K(Z_1,Z_2,Z_3)$ are exactly the same as the actions of $\sigma$ and
$\rho$ on $K(Y_1,Y_2,Y_3)$ in Case 3.

{\bf Case 12.} $\widetilde G=G_{12}$. We may again assume that
$a<t+\alpha$ and $\tau^{p^t}\in Z(\widetilde G)$. The subgroup
$H=\langle\sigma,\rho,\tau^{p^t}\rangle$ is abelian and has an order
$p^{a+b}$.

{\bf Subcase 12.a.} $b-\beta\leq t$. Then $c\geq\beta\geq b-t$,
i.e., $c-b+t\geq 0$. Define
$\rho_1=\sigma^{p^{a-\alpha}}\rho^{-1},\rho_2=\tau^{p^t}\sigma^{-p^{c+t-b}}\rho^{-p^{\beta+t-b}}$.
Then $H$ is isomorphic to the direct product
$\langle\sigma\rangle\times\langle\rho_1\rangle\times\langle\rho_2\rangle$,
where
$\sigma^{p^{a+t-\alpha}}=\rho_1^{p^\alpha}=\rho_2^{p^{b-t}}=1$. We
have $\rho=\rho_1^{-1}\sigma^{p^{a-\alpha}}$ and
$\tau^{p^t}=\rho_2\sigma^{p^{c+t-b}}\rho^{p^{\beta+t-b}}$.

{\bf Subcase 12.b.} $b-\beta>t$. We have
$(\tau^{p^{b-\beta}}\sigma^{-p^{c-\beta}}\rho^{-1})^{p^\beta}=1$.
Clearly, $\widetilde G'\cap
\langle\tau^{p^{b-\beta}}\sigma^{-p^{c-\beta}}\rho^{-1}\rangle=\{1\}$,
so we can apply Theorem \ref{t2.8} reducing the rationality problem
of $K(\widetilde G)$ to the rationality problem of $K(\widetilde
G_1)$ over $K$, where $\widetilde G_1\cong\widetilde
G/\langle\tau^{p^{b-\beta}}\sigma^{-p^{c-\beta}}\rho^{-1}\rangle$.
The group $\widetilde G_1$ is generated by elements $\sigma,\tau$
and $\rho$ such that
$\sigma^{p^a}=\rho^{p^\alpha},\tau^{p^{b-\beta}}=\sigma^{p^{c-\beta}}\rho,\rho^{p^t}=1,[\sigma,\tau]=\rho,
\rho - \text{central}$. Then the abelian subgroup
$H=\langle\sigma,\rho,\tau^{p^t}\rangle$ is of order
$p^{a+b-\beta}$.

{\bf Sub-subcase 12.b.a.} Let $b-c-t>0$ and $\alpha\geq\beta+a-c$.
Put $c-\beta+\alpha=a+u$, where $u\geq 0$. Define
$\rho_1=\sigma^{1+p^u}\tau^{-p^{b-c+u}}$ and $\rho_2=\tau^{p^t}$.
Then $H$ is isomorphic to the direct product
$\langle\rho_1\rangle\times\langle\rho_2\rangle$, where
$\rho_1^{p^a}=\rho_2^{p^{b-\beta}}=1$. Let $x$ be an integer such
that $(1+p^u)x\equiv 1\pmod{p^{a+t-\alpha}}$. Whence
$\sigma=(\sigma^{1+p^u})^x=\rho_1^x\rho_2^{xp^{b-c-t+u}}$ and
$\rho=\rho_1^{-xp^{c-\beta}}\rho_2^{p^{b-\beta-t}-xp^{b-t-\beta+u}}$.

{\bf Sub-subcase 12.b.b.} Let $b-c-t>0$ and $\alpha<\beta+a-c$. Put
$\beta+a-c=\alpha+v$, where $v>0$. Define
$\rho_1=\sigma^{1+p^v}\tau^{-p^{b-c}}$ and $\rho_2=\tau^{p^t}$. Then
$H$ is isomorphic to the direct product
$\langle\rho_1\rangle\times\langle\rho_2\rangle$, where
$\rho_1^{p^{a-v}}=\rho_2^{p^{b-\beta+v}}=1$. Let $x$ be an integer
such that $(1+p^v)x\equiv 1\pmod{p^{a+t-\alpha}}$. Whence
$\sigma=(\sigma^{1+p^v})^x=\rho_1^x\rho_2^{xp^{b-c-t}}$ and
$\rho=\rho_1^{-xp^{c-\beta}}\rho_2^{p^{b-\beta-t}-xp^{b-t-\beta}}$.

{\bf Sub-subcase 12.b.c.} Let $b-c-t\leq 0$ and
$y=b+\alpha-\beta-t>a$. Put $y=a+u$ and $c+t=b+v$, where $u>0,v\geq
0$. Note that $c+\alpha-\beta\geq b-t+\alpha-\beta=y>a$. Define
$\rho_2=\tau^{p^t}\sigma^{-p^v}$ and
$\rho_1=\sigma^{1+p^{c+\alpha-\beta-a}}\tau^{-p^{t+u}}=\sigma\rho_2^{-p^u}$.
Then $H$ is isomorphic to the direct product
$\langle\rho_1\rangle\times\langle\rho_2\rangle$, where
$\rho_1^{p^a}=\rho_2^{p^{b-\beta}}=1$. Note that
$\sigma=\rho_1\rho_2^{p^u},
\tau^{p^t}=\rho_1^{p^v}\rho_2^{1+p^{u+v}}$ and
$\rho=\rho_2^{p^{b-\beta-t}}$.

{\bf Sub-subcase 12.b.d.} Let $b-c-t\leq 0$ and
$y=b+\alpha-\beta-t\leq a$. Put $c+t=b+v$, where $v\geq 0$. Define
$\rho_1=\tau^{p^t}\sigma^{-p^v}\sigma^{-p^{a-y}}$. Then $H$ is
isomorphic to the direct product
$\langle\sigma\rangle\times\langle\rho_1\rangle$, where
$\sigma^{p^{a+t-\alpha}}=\rho_1^{p^y}=1$. Note that
$\tau^{p^t}=\rho_1\sigma^{p^v+p^{a-y}}$ and
$\rho=\rho_1^{p^{b-\beta-t}}\sigma^{p^{a-\alpha}}$.

\medskip
{\bf Remark.} For Sub-subcases 12.b.a and 12.b.b we have
$\exp(\widetilde G)\geq b>c+t\geq r+t$, so we can also apply Theorem
\ref{t1.5}.
\medskip

 {\bf Step II.} Assume that $\widetilde G$ has the following
general presentation:
\begin{equation*}
\widetilde
G=\langle\sigma,\tau,\rho:\sigma^{p^a}=\rho^{sp^\alpha},\tau^{p^b}=\sigma^{mp^c}\rho^{p^\beta},\rho^{p^t}=1,\tau^{-1}\sigma\tau=\sigma^k\rho,
\rho - \text{central}\rangle,
\end{equation*}
where $s,m$ are positive integers, $1\leq s<p^t,0 \leq m<p^a$,
$\gcd(sm,p)=1,0\leq\alpha,\beta\leq t$ and $k=\varepsilon+p^r$.

Here we have again $16$ cases, which correspond to those in Step I.
Since they all can be treated in an unified way, we will consider
only Case 3.

Let $c=\beta=t$, i.e., $\tau^{p^b}=1$.

We have the relations
$(\sigma^{p^{a-\alpha}}\rho^{-s})^{p^\alpha}=1$ and
$\tau^{-1}\sigma^{p^{a-\alpha}}\tau=\sigma^{p^{a-\alpha}}\rho^{p^{a-\alpha}}$.
Clearly, $\widetilde
G'\cap\langle\sigma^{p^{a-\alpha}}\rho^{-s}\rangle=\{1\}$. If
$a-\alpha\geq t$ then $\sigma^{p^{a-\alpha}}\in Z(\widetilde G)$.
Theorem \ref{t2.8} then implies that we can reduce the rationality
problem of $K(\widetilde G)$ to $K(\widetilde
G/\langle\sigma^{p^{a-\alpha}}\rho^{-s}\rangle)$ over $K$, where
$\widetilde G_1=\widetilde
G/\langle\sigma^{p^{a-\alpha}}\rho^{-s}\rangle$ is a metacyclic
$p$-group. Indeed, if $\widetilde G_1$ is generated by elements
$\sigma,\tau,\rho$ such that $\sigma^{p^{a-\alpha}}=\rho^s$ and
$\tau^{-1}\sigma\tau=\sigma^k\rho$, then
$\tau^{-1}\sigma\tau=\sigma^{k+np^{a-\alpha}}$ for an integer $n$
with $\rho=(\rho^s)^n$, i.e., $ns\equiv 1\pmod{p^t}$.

Therefore, we will assume henceforth that $a-\alpha<t$. We are going
to show now that $\tau^{p^t}\in Z(\widetilde G)$. Indeed, we have
\begin{equation*}
\tau^{-p^t}\sigma\tau^{p^t}=\sigma^{(1+p^r)^{p^t}}=\sigma\cdot\sigma^{Ap^{r+t}}=\sigma\cdot\rho^{Ap^{\alpha+r+t-a}}=\sigma,
\end{equation*}
where $A$ is some integer and $\alpha+r+t-a>t$ since $a<\alpha+t\leq
\alpha+r$.

In this way, we obtain that the subgroup
$H=\langle\sigma,\rho,\tau^{p^t}\rangle$ is abelian and has an order
$p^{a+b}$. Put $\rho_1=\sigma^{p^{a-\alpha}}\rho^{-s}$ and
$\rho_2=\tau^{p^t}$. Then
$H\cong\langle\sigma\rangle\times\langle\rho_1\rangle\times\langle\rho_2\rangle$,
where
$\sigma^{p^{a+t-\alpha}}=\rho_1^{p^\alpha}=\rho_2^{p^{b-t}}=1$. Note
that $\rho^s=\sigma^{p^{a-\alpha}}\rho_1^{-1}$, so
$\rho=\sigma^{np^{a-\alpha}}\rho_1^{-n}$ with  $ns\equiv
1\pmod{p^t}$.

Define $X_1,X_2,X_3\in V^*$ by
\begin{equation*}
X_1=\sum_{i,j}x(\rho_1^i\rho_2^j),~
X_2=\sum_{i,j}x(\sigma^i\rho_1^j),~
X_3=\sum_{i,j}x(\sigma^i\rho_2^j).
\end{equation*}
Note that $\sigma\cdot X_2=X_2,\sigma\cdot X_3=X_3,\rho_1\cdot
X_1=X_1,\rho_1\cdot X_2=X_2,\rho_2\cdot X_1=X_1$ and $\rho_2\cdot
X_3=X_3$.

Let $\zeta_1\in K$ be a primitive $p^{a+t-\alpha}$-th root of unity.
Define $\zeta_3=\zeta_1^{p^{a+2t-\alpha-b}}$, a primitive
$p^{b-t}$-th root of unity; and $\xi=\zeta_1^{p^{a-\alpha}}$, a
primitive $p^t$-th root of unity. Since $a-\alpha<t\leq r$, we can
write $r=a-\alpha+r_1$ for some $r_1>0$. Define
$\zeta_2=\xi^{(1+sp^{r_1})p^{t-\alpha}}$, a primitive $p^\alpha$-th
root of unity (recall that $ns\equiv 1\pmod{p^t}$).

Define $Y_1,Y_2,Y_3\in V^*$ by
\begin{equation*}
Y_1=\sum_{i=0}^{p^{a+t-\alpha}-1}\zeta_1^{-i}\sigma^i\cdot X_1,~
Y_2=\sum_{i=0}^{p^\alpha-1}\zeta_2^{-i}\rho_1^i\cdot X_3,~
Y_3=\sum_{i=0}^{p^{b-t}-1}\zeta_3^{-i}\rho_2^i\cdot X_2.
\end{equation*}

It follows that {\allowdisplaybreaks \begin{align*}
\sigma\ :\ &Y_1\mapsto\zeta_1 Y_1,~ Y_2\mapsto Y_2,~ Y_3\mapsto Y_3,\\
\rho_1\ :\ &Y_1\mapsto Y_1,~ Y_2\mapsto\zeta_2 Y_2,~ Y_3\mapsto Y_3,\\
\rho_2\ :\ &Y_1\mapsto Y_1,~ Y_2\mapsto Y_2,~ Y_3\mapsto\zeta_3 Y_3,\\
\rho\ :\ &Y_1\mapsto\xi^n Y_1,~ Y_2\mapsto\zeta_2^{-n} Y_2,~
Y_3\mapsto Y_3.
\end{align*}}
Thus $K\cdot Y_1+K\cdot Y_2+K\cdot Y_3$ is a representation space of
the subgroup $H$.

Define $x_i=\tau^i\cdot Y_1,y_i=\tau^i\cdot Y_2,z_i=\tau^i\cdot Y_3$
for $0\leq i\leq p^t-1$. We have now {\allowdisplaybreaks
\begin{align*}
\sigma\ :\ &x_i\mapsto \zeta_1^{k^i}\xi^{ni} x_i,~ y_i\mapsto\zeta_2^{-ni}y_i,~ z_i\mapsto z_i,\\
\tau\ :\ &x_0\mapsto x_1\mapsto\cdots\mapsto x_{p^t-1}\mapsto x_0,\\
&y_0\mapsto y_1\mapsto\cdots\mapsto y_{p^t-1}\mapsto y_0,\\
&z_0\mapsto z_1\mapsto\cdots\mapsto z_{p^t-1}\mapsto\zeta_3 z_0,\\
\rho\ :\ &x_i\mapsto\xi^n x_i,~ y_i\mapsto\zeta_2^{-n} y_i,~
z_i\mapsto z_i.
\end{align*}}
for $0\leq i\leq p^t-1$. We find that $Y=(\bigoplus_{0\leq i\leq
p^t-1}K\cdot x_i)\oplus(\bigoplus_{0\leq i\leq p^t-1}K\cdot
y_i)\oplus(\bigoplus_{0\leq i\leq p^t-1}K\cdot z_i)$ is a faithful
$\widetilde G$-subspace of $V^*$. Thus, by Theorem \ref{t2.1}, it
suffices to show that $K(x_i,y_i,z_i:0\leq i\leq p^t-1)^{\widetilde
G}$ is rational over $K$.

For $1\leq i\leq p^t-1$, define $U_i=x_i/x_{i-1}$ and
$V_i=y_i/y_{i-1}$. Thus $K(x_i,y_i,z_i:0\leq i\leq
p^t-1)=K(x_0,y_0,U_i,V_i,z_j:1\leq i\leq p^t-1,0\leq j\leq p^t-1)$
and for every $g\in\widetilde G$
\begin{equation*}
g\cdot x_0\in K(U_i,V_i,z_j)\cdot x_0,~ g\cdot y_0\in
K(U_i,V_i,z_j)\cdot y_0,
\end{equation*}
while the subfield $K(U_i,V_i,z_j:1\leq i\leq p^t-1,0\leq j\leq
p^t-1)$ is invariant by the action of $\widetilde G$, i.e.,
{\allowdisplaybreaks
\begin{align*}
\sigma\ :\ &U_i\mapsto \zeta_1^{k^i-k^{i-1}}\xi^n U_i,~ V_i\mapsto\zeta_2^{-n} V_i,~ z_j\mapsto z_j,\\
\tau\ :\ &U_1\mapsto U_2\mapsto\cdots\mapsto U_{p^t-1}\mapsto (U_1\cdots U_{p^t-1})^{-1},\\
&V_1\mapsto V_2\mapsto\cdots\mapsto V_{p^t-1}\mapsto (V_1\cdots V_{p^t-1})^{-1},\\
&z_0\mapsto z_1\mapsto\cdots\mapsto z_{p^t-1}\mapsto\zeta_3 z_0,\\
\rho\ :\ &U_i\mapsto U_i,~ V_i\mapsto V_i,~ z_j\mapsto z_j.
\end{align*}}
for $1\leq i\leq p^t-1$ and $0\leq j\leq p^t-1$. From Theorem
\ref{t2.2} it follows that if $K(U_i,V_i,z_j:1\leq i\leq p^t-1,0\leq
j\leq p^t-1)^{\widetilde G}$ is rational over $K$, so is
$K(x_i,y_i,z_i:0\leq i\leq p^t-1)^{\widetilde G}$ over $K$.

Since $\rho$ acts trivially on $K(U_i,V_i,z_j)$, we find that
$K(U_i,V_i,z_j)^{\widetilde
G}=K(U_i,V_i,z_j)^{\langle\sigma,\tau\rangle}$.

Recall that $r_1=r-a+\alpha>0$. Therefore,
$\zeta_1^{k-1}=\zeta_1^{p^{a-\alpha+r_1}}=\xi^{p^{r_1}}$ and also
$\zeta_1^{k^{i-1}(k-1)}=\xi^{p^{r_1}}$ for all $i$.

Define $v_i=U_i^{p^{t-\alpha}}V_i$. Since
$\xi^{(n+p^{r_1})p^{t-\alpha}}=\xi^{n(1+sp^{r_1})p^{t-\alpha}}=\zeta_2^n$,
we have {\allowdisplaybreaks
\begin{align*}
\sigma\ :\ &U_i\mapsto \xi^{n+p^{r_1}} U_i,~ v_i\mapsto v_i,~ z_j\mapsto z_j,\\
\tau\ :\ &U_1\mapsto U_2\mapsto\cdots\mapsto U_{p^t-1}\mapsto (U_1\cdots U_{p^t-1})^{-1},\\
&v_1\mapsto v_2\mapsto\cdots\mapsto v_{p^t-1}\mapsto (v_1\cdots v_{p^t-1})^{-1},\\
&z_0\mapsto z_1\mapsto\cdots\mapsto z_{p^t-1}\mapsto\zeta_3 z_0.
\end{align*}}
for $1\leq i\leq p^t-1$ and $0\leq j\leq p^t-1$. The actions of
$\sigma$ and $\tau$ now are very similar to the actions in Case 3,
Step I. Apply the same proof.

\end{document}